\documentclass[11pt,reqno]{amsart}
\usepackage{fullpage,times,graphicx,amssymb,amsmath,amsfonts,psfrag,xcolor,pdfsync,stmaryrd,enumitem}
\usepackage[colorlinks,citecolor=bbluegray,linkcolor=ddarkbrown,urlcolor=blue,breaklinks]{hyperref}
\newtheorem{theorem}{Theorem}[section]
\newtheorem{proposition}[theorem]{Proposition}
\newtheorem{definition}[theorem]{Definition}
\newtheorem{lemma}[theorem]{Lemma}
\newtheorem{corollary}[theorem]{Corollary}
\renewenvironment{proof}{\textbf{Proof.}}{\QED\bigskip}

\usepackage{algorithm,algorithmicx,algpseudocode}

\usepackage[square]{natbib}

\definecolor{ddarkbrown}{rgb}{0.5,0.2,0.05} \definecolor{bbluegray}{rgb}{0.05,0,0.5}

\newcommand{\BEAS}{\begin{eqnarray*}}
\newcommand{\EEAS}{\end{eqnarray*}}
\newcommand{\BEA}{\begin{eqnarray}}
\newcommand{\EEA}{\end{eqnarray}}
\newcommand{\BEQ}{\begin{equation}}
\newcommand{\EEQ}{\end{equation}}
\newcommand{\BIT}{\begin{itemize}}
\newcommand{\EIT}{\end{itemize}}
\newcommand{\BNUM}{\begin{enumerate}}
\newcommand{\ENUM}{\end{enumerate}}

\newcommand{\BA}{\begin{array}}
\newcommand{\EA}{\end{array}}

\newcommand{\reals}{{\mathbb R}}

\newcommand{\Supp}{\mbox{\textrm{Supp}}}

\newcommand{\diam}{\mathop{\bf diam}}

\newcommand{\Rank}{\mathop{\bf Rank}}

\newcommand{\Card}{\mathop{\bf Card}}

\newcommand{\idm}{\mathbf{I}}

\newcommand{\QED}{~~\rule[-1pt]{6pt}{6pt}}

\newcommand{\intr}{\mathop{\bf int}}

\renewcommand\Im{\operatorname{Im}}

\usepackage{auto-pst-pdf} 

\newcommand{\calE}{E}

\newcommand{\Null}{\mathrm{Null}}
\let \oldsection \section
\renewcommand{\section}{\vspace{3ex plus 1ex}\oldsection}

\title{Computational Complexity versus Statistical Performance \\on Sparse Recovery Problems}

\author{Vincent Roulet} 
\address{DI, \'Ecole Normale Sup\'erieure, Paris, France.\vskip 0ex 
	INRIA Sierra Team.}
\email{vincent.roulet@inria.fr}

\author{Nicolas Boumal} 
\address{Mathematics Department, \vskip 0ex
	Princeton university, Princeton NJ 08544, USA.}
\email{nboumal@math.princeton.edu}

\author{Alexandre d'Aspremont} 
\address{CNRS \& DI, \'Ecole Normale Sup\'erieure, Paris, France.\vskip 0ex}
\email{aspremon@di.ens.fr}

\date{\today}

\begin{document}

\maketitle
\begin{abstract} 
We show that several classical quantities controlling compressed sensing performance directly match classical parameters controlling algorithmic complexity. We first describe linearly convergent restart schemes on first-order methods solving a broad range of compressed sensing problems, where sharpness at the optimum controls convergence speed. We show that for sparse recovery problems, this sharpness can be written as a condition number, given by the ratio between true signal sparsity and the largest signal size that can be recovered by the observation matrix. In a similar vein, Renegar's condition number is a data-driven complexity measure for convex programs, generalizing classical condition numbers for linear systems. We show that for a broad class of compressed sensing problems, the worst case value of this algorithmic complexity measure taken over all signals matches the restricted singular value of the observation matrix which controls robust recovery performance. Overall, this means in both cases that, in compressed sensing problems, a single parameter directly controls both computational complexity and recovery performance. Numerical experiments illustrate these points using several classical algorithms.
\end{abstract}

\section*{Introduction}
Sparse recovery problems have received a lot of attention from various perspectives. On one side, an extensive literature explores the limits of recovery performance. On the other side, a long list of algorithms now solve these problems very efficiently. Early on, it was noticed empirically by e.g.~\cite{Dono06}, that recovery problems which are easier to solve from a statistical point of view (i.e., where more samples are available), are also easier to solve numerically. Here, we show that these two aspects are indeed intimately related. 

Recovery problems consist in retrieving a signal $x^*$, lying in some Euclidean space $E$, given linear observations. If the signal is ``sparse", namely if it can be efficiently compressed, a common approach is to minimize the corresponding sparsity inducing norm $\|\cdot\|$ (e.g. the $\ell_1$ norm in classical sparse recovery). The exact sparse recovery problem then reads
\BEQ\label{eq:recov}
\BA{ll}
\mbox{minimize} & \|x\|\\
\mbox{subject to} & A(x)=b,
\EA
\EEQ
in the variable $x\in E$, where $A$ is a linear operator on $E$ and $b = A(x^*)$ is the vector of observations. If the observations are affected by noise a robust version of this problem is written as
\BEQ
\BA{ll}
\mbox{minimize} & \|x\|\\
\mbox{subject to} & \|A(x)-b\|_2 \leq \epsilon,
\EA
\EEQ
in the variable $x\in E$, where $\|\cdot\|_2$ is the Euclidean norm and $\epsilon >0$ is a tolerance to noise. In penalized form, this is
\BEQ
\BA{ll}
\mbox{minimize} & \|x\| + \lambda\|A(x)-b\|^2_2
\EA
\EEQ
in the variable $x\in E$ where $\lambda > 0$ is a penalization parameter. This last problem is known as the LASSO \citep{Tibs96} in the $\ell_1$ case.

When $x^*$ has no more than $s$ non zero values, \cite{Dono05} and \cite{Cand06} have shown that, for certain linear operators $A$, $O(s\log p)$ observations suffice for stable recovery of $x^*$ by solving the exact formulation~\eqref{eq:recov} using the $\ell_1$ norm (a linear program), where $p$ is the dimension of the space $E$. These results have then been generalized to many other recovery problems with various assumptions on signal structure (e.g., where $x$ is a block-sparse vector, a low-rank matrix, etc.) and corresponding convex relaxations were developed in those cases (see e.g. \cite{Chan12} and references therein). Recovery performance is often measured in terms of the number of samples required to guarantee exact or robust recovery given a level of noise. 

On the computational side, many algorithms were developed to solve these problems at scale. Besides specialized methods such as LARS \citep{Efro04}, FISTA \citep{Beck09} and NESTA \citep{Beck11}, solvers use accelerated gradient methods to solve robust recovery problems, with efficient and flexible implementations covering a wide range of compressed sensing instances developed by e.g. \cite{Beck12}. Recently, linear convergence results have been obtained for the LASSO \citep{Agar11,Yen14,Zhou15} using variants of the classical strong convexity assumption, while \citep{Zhou17} studied error bounds for a much broader class of structured optimization problems including sparse recovery and matrix completion. Some restart schemes have also been developed in e.g. \citep{Odon15,Su14,Gise14} while \cite{Ferc16} showed that generic restart schemes can offer linear convergence given a rough estimate of the behavior of the function around its minimizers.

As mentioned above, \cite{Dono06} was one of the first reference to connect statistical and computational performance in this case, showing empirically that recovery problems which are easier to solve from a statistical point of view (i.e., where more samples are available), are also easier to solve numerically (using homotopy methods). More recently, \cite{Chan13,Amel14} studied computational and statistical tradeoffs for increasingly tight convex relaxations of shrinkage estimators. They show that recovery performance is directly linked to the Gaussian squared-complexity of the tangent cone with respect to the constraint set and study the complexity of several convex relaxations. In \citep{Chan13,Amel14} however, the structure of the convex relaxation is varying and affecting both complexity and recovery performance, while in \citep{Dono06} and in what follows, the structure of the relaxation is fixed, but the data (i.e. the observation matrix $A$) varies.

Here, as a first step, we study the exact recovery case and show that the null space property introduced by \cite{Cohe06} can be seen as a measure of sharpness on the optimum of the sparse recovery problem. On one hand this allows us to develop linearly convergent restart schemes whose rate depends on this sharpness. On the other hand we recall how the null space property is linked to the recovery threshold of the sensing operator $A$ for random designs, thus producing a clear link between statistical and computational performance.

We then analyze the underlying conic geometry of recovery problems. Robust recovery performance is controlled by a minimal conically restricted singular value. We recall Renegar's condition number and show how it affects the computational complexity of optimality certificates for exact recovery and the linear convergence rate of restart schemes. By observing that the minimal conically restricted singular value matches the worst case value of Renegar's condition number on sparse signals, we provide further evidence that a single quantity controls both computational and statistical aspects of recovery problems. Numerical experiments illustrate its impact on various classical algorithms for sparse recovery. 


The first two sections focus on the $\ell_1$ case for simplicity. We generalize our results to non-overlapping group norms and the nuclear norm in a third section.

\subsection*{Notations}
For a given integer $p\geq 1$, $\llbracket 1, p \rrbracket$ denotes the set of integers between $1$ and $p$. For a given subset $S\subset\llbracket1,p\rrbracket$, we denote $S^c= \llbracket 1, p \rrbracket \setminus S$ its complementary and $\Card(S)$ its cardinality. For a given vector $x\in \reals^p$, we denote $\Supp(x) = \{i \in \llbracket 1,p\rrbracket : x_i \neq 0\}$ the support of $x$, $\|x\|_0 = \Card(\Supp(x))$ its sparsity and $\|x\|_p$ its $p$-norm. For a given vector $x$ and integer subset $S \subset \llbracket 1, p \rrbracket$, $x_S \in \reals^p$ denotes the vector obtained by zeroing all coefficients of $x$ that are not in $S$.
For a given linear operator or matrix $A$, we denote $\Null(A)$ its null space, $\Im(A)$ its range, and $\|X\|_2$ its operator norm with respect to the Euclidean norm (for matrices this is the spectral norm). The identity operator is denoted $\idm$. In a linear topological space $E$ we denote $\intr(F)$ the interior of $F\subset E$. Finally for a given real $a$, we denote $\lceil a \rceil$ the smallest integer larger than or equal to $a$ and $\lfloor a \rfloor$ the largest integer smaller than or equal to $a$.

\section{Sharpness, Restart and Sparse Recovery Performance}\label{s:sharp}

In this section and the following one, we discuss sparse recovery problems using the $\ell_1$ norm. Given a matrix $A\in\reals^{n\times p}$ and observations $b=Ax^*$ on a signal $x^*\in\reals^p$, recovery is performed by solving the $\ell_1$ minimization program
\BEQ\label{eq:l1-recov}\tag{$\ell_1$ recovery}
\BA{ll}
\mbox{minimize} & \|x\|_1\\
\mbox{subject to} & Ax=b
\EA
\EEQ
in the variable $x\in\reals^p$. 

In what follows, we show that the Null Space Property condition (recalled below) can be seen as measure of sharpness for $\ell_1$-recovery of a sparse signal $x^*$, with
\BEQ\label{eq:sharpness}\tag{Sharp}
\|x\|_1-\|x^*\|_1 > \gamma \|x-x^*\|_1
\EEQ
for any $x \neq x^*$ such that $Ax=b$, and some $0 \leq \gamma <1$. This first ensures that $x^*$ is the unique minimizer of problem~\eqref{eq:l1-recov} but also has important computational implications. It allows us to produce linear convergent restart schemes whose rates depend on sharpness. By connecting null space property to recovery threshold for random observation matrices, we thus get a direct link between computational and statistical aspects of sparse recovery problems.

\subsection{Null space property \& sharpness for exact recovery}\label{s:sharp_nullspace}
Although the definition of null space property appeared in earlier work \citep{Dono01,Feue03} the terminology of restricted null space is due to \cite{Cohe06}. The following definition differs slightly from the original one in order to relate it to intrinsic geometric properties of the problem in Section~\ref{s:conic_linear}.

\begin{definition}{\bf (Null Space Property) }\label{def:nsp_prop} The matrix $A$ satisfies the Null Space Property (NSP) \emph{on support $S\subset \llbracket 1,p \rrbracket$} with constant $\alpha \geq 1$ if for any $z\in \Null(A)\setminus\{0\}$,
	\BEQ\label{def:nsp}\tag{NSP}
	\alpha \|z_S\|_1 < \|z_{S^c}\|_1.
	\EEQ
	The matrix $A$ satisfies the Null Space Property \emph{at order $s$} with constant $\alpha \geq 1$ if it satisfies it on every support $S$ of cardinality at most $s$.
\end{definition}

The Null Space Property is a necessary and sufficient condition for the convex program \eqref{eq:l1-recov} to recover all signals up to some sparsity threshold. Necessity will follow from results recalled in Section~\ref{s:recov}. We detail sufficiency of \eqref{def:nsp} using sharpness in the following proposition.

\begin{proposition}\label{prop:recov}
	Given a coding matrix $A\in\reals^{n \times p}$ satisfying \eqref{def:nsp} at order $s$ with constant $\alpha \geq 1$, if the original signal $x^*$ is $s$-sparse, then for any $x\in\reals^p$ satisfying $Ax=b$, $x\neq x^*$, we have  
	\BEQ
	\|x\|_1-\|x^*\|_1 > \frac{\alpha-1}{\alpha+1} \|x-x^*\|_1.
	\EEQ
	This implies signal recovery, i.e. optimality of $x^*$ for \eqref{eq:l1-recov}, and the sharpness bound \eqref{eq:sharpness} with $\gamma = \frac{\alpha-1}{\alpha+1}$.
\end{proposition}
\begin{proof}
	The proof follows the one in \cite[Theorem 4.4]{Cohe06}. Let $S=\mathrm{supp}(x^*)$, with $\Card(S) \leq s$, and let $x\neq x^*$ such that $Ax=b$, so $z=x-x^*\neq 0$ satisfies $Az=0$. Then
	\begin{align}
	\|x\|_1 & =  \|x^*_{S}+z_{S}\|_1 + \|z_{S^c}\|_1 \\ 
	& \geq  \|x^*_{S}\|_1 - \|z_{S}\|_1 + \|z_{S^c}\|_1 \\ 
	& =  \|x^*\|_1 + \|z\|_1 - 2 \|z_{S}\|_1 . 
	\end{align}
	Now as $A$ satisfies \eqref{def:nsp} on support $S$, 
	\BEQ
	\|z\|_1 = \|z_S\|_1 + \|z_{S^c}\|_1 > (1+\alpha)\|z_{S}\|_1
	\EEQ 
	hence
	\begin{align}
	\|x\|_1 - \|x^*\|_1 > \frac{\alpha-1}{\alpha+1} \|z\|_1 = \frac{\alpha-1}{\alpha+1} \|x-x^*\|_1 .
	\end{align}
	As $\alpha \geq 1$, this implies that $x^*$ is the solution of program~\eqref{eq:l1-recov} and the corresponding sharpness bound.
\end{proof}

Sharpness is a crucial property for optimization problems that can be exploited to accelerate the performance of classical optimization algorithms \citep{Nemi85,Roul17}. Before that we remark that it is in fact equivalent to \eqref{def:nsp} at order $s$.
\begin{proposition}\label{prop:sharp_to_nsp}
	Given a matrix $A\in \reals^{n\times p}$ such that problem~\eqref{eq:l1-recov} is sharp on every $s$-sparse signal $x^*$, i.e. there exists $ 0 \leq \gamma <1$ such that
	\BEQ
	\|x\|_1-\|x^*\|_1 > \gamma \|x-x^*\|_1,
	\EEQ
	for any $x\neq x^*$ such that $Ax = Ax^*$. Then, $A$ satisfies \eqref{def:nsp} at order $s$ with constant $\alpha = \frac{1 + \gamma}{1-\gamma} \geq 1$.
\end{proposition}
\begin{proof}
	Let $S\subset \llbracket 1,p\rrbracket$ with $\Card(S) \leq s$ and $z \in \Null(A)$, $z\neq0$,  such that $Az_S = -Az_{S^c}$ and $z_S \neq -z_{S^c}$. Using sharpness of problem~\eqref{eq:l1-recov} with $x^* = z_S$, and $x = -z_{S^c}$, we get
	\BEQ
	\|z_{S^c}\|_1-\|z_S\|_1 > \gamma \|z\|_1 = \gamma \|z_S\|_1 + \gamma\|z_{S^c}\|_1.
	\EEQ
	Rearranging terms and using $\gamma<1$, this reads 
	\BEQ
	\|z_{S^c}\|_1 >\frac{1 + \gamma}{1-\gamma} \|z_S\|_1,
	\EEQ
	which is \eqref{def:nsp} on support $S$ with the announced constant. As $S$ was taken arbitrarily, this means \eqref{def:nsp} holds  at order $s$.
\end{proof}

\subsection{Restarting first-order methods} \label{s:restart}
In this section, we seek to solve the recovery problem \eqref{eq:l1-recov} and exploit the sharpness bound~\eqref{eq:sharpness}. The NESTA algorithm \citep{Beck11} uses the smoothing argument of \cite{Nest03} to solve \eqref{eq:l1-recov}. In practice, this means using the optimal algorithm of \cite{Nest83} to minimize 
\BEQ
f_\epsilon(x) \triangleq \sup_{\|u\|_\infty\leq 1} \left\{u^Tx - \frac{\epsilon}{2p} \|u\|_2^2\right\}
\EEQ
for some $\epsilon>0$, which approximates the $\ell_1$ norm uniformly up to $\epsilon/2$. This is the classical Huber function, which has a Lipschitz continuous gradient with constant equal to $p/\epsilon$. Overall given an accuracy $\epsilon$ and a starting point $x_0$ this method outputs after $t$ iterations a point $x = \mathcal{A}(x_0,\epsilon,t)$ such that 
\BEQ\label{eq:conv_rate}
\|x\|_1-\|\hat x\|_1 \leq \frac{2p \|x_0-\hat{x}\|_2^2}{\epsilon t^2} + \frac{\epsilon}{2},
\EEQ
for any $\hat x$ solution of problem~\eqref{eq:l1-recov}.
Now if the sharpness bound is satisfied, restarting this method, as described in the \eqref{eq:restart} scheme presented below, accelerates its convergence. 

\begin{algorithm}[h]
	\caption{Restart Scheme \textbf{(Restart)}}
	\label{algo:restart}
	\begin{algorithmic}
		\Require Initial point $y_0\in\reals^p$, initial gap $\epsilon_0\geq \|y_0\|_1-\|\hat x\|_1$, decreasing factor $\rho$, restart clock $t$
		\State For $k=1\ldots,K$ compute
		\BEQ\label{eq:restart}\tag{Restart}
		\epsilon_k = \rho \epsilon_{k-1}, \qquad y_k=\mathcal{A}(y_{k-1}, \epsilon_k, t)
		\EEQ
		\Ensure A point $\hat y = y_K$ approximately solving~\eqref{eq:l1-recov}.
	\end{algorithmic}
\end{algorithm}

\subsubsection{Optimal restart scheme}
We begin by analyzing an optimal restart scheme assuming the sharpness constant is known. We use a non-integer clock to highlight its dependency to the sharpness. Naturally clock and number of restarts must be integer but this does not affect much bounds as detailed in Appendix~\ref{app:practical_optimal}.
The next proposition shows that algorithm $\mathcal{A}$ needs a constant number of iterations to decrease the gap by a constant factor, which means restart leads to linear convergence. 

\begin{proposition}\label{prop:complex-l1}
Given a coding matrix $A\in\reals^{n \times p}$ and a signal $x^* \in \reals^p$ such that the sharpness bound~\eqref{eq:sharpness} is satisfied with $\gamma >0$, i.e.
\BEQ
\|x\|_1-\|x^*\|_1 > \gamma \|x-x^*\|_1,
\EEQ
for any $x\neq x^*$ such that $Ax = Ax^*$, running the \eqref{eq:restart} scheme with $t \geq \frac{2 \sqrt{p}}{\gamma \rho}$
ensures 
\BEQ\label{eq:goal_restart}
\|y_k\|_1 - \|x^*\|_1 \leq \epsilon_k,
\EEQ
at each iteration, with $x^*$ the unique solution of problem~\eqref{eq:l1-recov}.
Using optimal parameters
\BEQ\label{eq:optimal_param}
\rho^* = e^{-1} \qquad and \qquad t^* = \frac{2 e\sqrt{p}}{\gamma},
\EEQ
we get a point $\hat y$ such that 
\BEQ\label{eq:complexity}
\|\hat y\|_1-\|x^*\|_1 \leq \exp\left( -\frac{\gamma}{2 \sqrt{p}} e N  \right) \epsilon_0.
\EEQ
after running a total of $N$ inner iterations of Algorithm~\ref{eq:restart} with $t = t^*$ (hence $N/t$ restarts).
\end{proposition}
\begin{proof}
By the choice of $\epsilon_0$, \eqref{eq:goal_restart} is satisfied for $k=0$. Assuming it holds at iteration $k$, combining sharpness bound~\eqref{eq:sharpness} and complexity bound~\eqref{eq:conv_rate} leads to, for $x = \mathcal{A}(y_{k-1},\epsilon_k,t)$,
\begin{align*}
\|x\|_1-\|x^*\|_1 & \leq \frac{2p (\|y_{k-1}\|_1 - \|x^*\|_1)^2}{\gamma^2\epsilon_k t^2} + \frac{\epsilon_k}{2} \\
& \leq \frac{4p }{\rho^2\gamma^2 t^2} \frac{\epsilon_k}{2}+ \frac{\epsilon_k}{2}.
\end{align*}
Therefore after 
$t \geq \frac{2\sqrt{p}}{\gamma\rho} $
iterations, the method has achieved the decreased accuracy $\epsilon_k$ which proves \eqref{eq:goal_restart}. The overall complexity after a total of $N$ inner iterations, hence $ N/t$ restarts, is then
\BEQ
\|\hat y\|_1-\|x^*\|_1 \leq \rho^{N/t} \epsilon_0.
\EEQ
If $\gamma$ is known, using exactly $\frac{2 \sqrt{p}}{\gamma\rho}$ inner iterations at each restart leads to
\BEQ
\|\hat y\|_1-\|x^*\|_1 \leq \exp\left( \frac{\gamma}{2\sqrt{p}}N\rho \log \rho \right) \epsilon_0.
\EEQ
Optimizing in $\rho$ yields $\rho^* = e^{-1}$, and with $t^*$ inner iterations the complexity bound \eqref{eq:complexity} follows. 
\end{proof}
	
To run NESTA, $A^TA$ is assumed to be an orthogonal projector (w.l.o.g. at the cost of computing an SVD) such that the projection on the feasible set is easy. \cite{Beck11} already studied restart schemes that they called ``acceleration with continuation". However their restart criterion depends on the relative variation of objective values, not on the number of iterates, and no linear convergence was proven. We further note that linear convergence of restart schemes requires an assumption of the form
\BEQ\label{eq:Loja}
f(x) - f^* \geq \gamma d(x,X^*)^\nu,
\EEQ
where $d(x,X^*)$ is the distance (in any norm) from $x$ to the set of minimizers of the objective function $f$ (here $f(x) =\|x\|_1$). This type of bound is known as {\L}ojasiewicz's inequality, studied for example in \cite{Bolt07} for non-smooth convex functions. Here \eqref{def:nsp} ensures that the set of minimizers is reduced to a singleton, the original signal.

\subsubsection{Practical restart scheme}\label{s:pract_rest}
Several parameters are needed to run the optimal scheme above. The optimal decreasing factor is independent of the data. The initial gap $\epsilon_0$ can be taken as $\|y_0\|_1$ for $Ay_0 =b$. The sharpness constant $\gamma$ is for its part mostly unknown such that we cannot choose the number $t^*$ of inner iterations a priori. However, given a budget of iterations $N$ (the total number of iterations in the optimization algorithm, across restarts), a log scale grid search can be performed on the optimal restart clock to get nearly optimal rates as detailed in the following corollary (contrary to the general results in \citep{Roul17}, the sharpness exponent $\nu$ in~\eqref{eq:Loja} is equal to one here, simplifying the parameter search). 

\begin{corollary}\label{cor:grid_search}
Given a coding matrix $A\in\reals^{n \times p}$, a signal $x^* \in \reals^p$ such that the sharpness bound~\eqref{eq:sharpness} is satisfied with $\gamma >0$, a budget of $N$ iterations, run the following schemes from an initial point $y_0$
\BEQ
\text{\eqref{eq:restart} with } \quad t = h^j, \qquad  j= 1,\ldots, \lfloor\log_h N\rfloor 
\EEQ
with $h$ the grid search precision. Stop restart iteration when the total number of iterations has exceeded the budget $N$.
Then, provided that $N \geq h t^*$, where $t^*$ is defined in \eqref{eq:optimal_param}, at least one of these restart schemes achieves a precision given by
\BEQ\label{eq:complexity_grid_search}
\|\hat y\|_1-\|x^*\|_1 \leq \exp\left( -\frac{\gamma}{2 h \sqrt{p}} e N  \right) \epsilon_0.
\EEQ
Overall running the logarithmic grid search has a complexity $\log_h N$ times higher than running $N$ iterations in the optimal scheme.
\end{corollary}
\begin{proof}
All schemes stop after at most $N +h^j \leq 2N$ iterations. As we assumed  $N \geq h t^*$, $j = \lceil \log_h t^* \rceil \leq \log_h N$ and \eqref{eq:restart} has been run with $t = h^j$. Proposition~\ref{prop:complex-l1} ensures, since $t\geq t^*$, that the output of this scheme achieves after $N'\geq N$ total iterations a precision 
\BEQ
\|\hat y\|_1 - \|\hat x\|_1 \leq e^{-N'/t} \epsilon_0 \leq  e^{-N/t} \epsilon_0
\EEQ
and as $t \leq h t^*$
\BEQ
\|\hat y\|_1 - \|\hat x\|_1 \leq e^{-N/ (h t^*)} \epsilon_0
\EEQ
which gives the result. Finally the logarithmic grid search costs $\log_h N$ to get this approximative optimal bound. 
\end{proof}

Sharpness therefore controls linear convergence of simple restart schemes to solve \eqref{eq:l1-recov}. We now turn back to \eqref{def:nsp} estimates and connect them to recovery thresholds of the sampling matrix. This will give us a direct link between computational complexity and recovery performance for exact recovery problems.

\subsection{Recovery threshold}\label{s:oversampling}
If \eqref{def:nsp} is satisfied at a given order $s$ it holds also for any $s'\leq s$. However, the constant, and therefore the speed of convergence, may change. Here we show that this constant actually depends on the ratio between the maximal order at which $A$ satisfies \eqref{def:nsp} and the sparsity of the signal that we seek to recover.

To this end, we give a more concrete geometric meaning to the constant $\alpha$ in~\eqref{def:nsp}, connecting it with the diameter of a section of the $\ell_1$ ball by the null space of the matrix $A$ (see e.g. \cite{Kash07} for more details).

\begin{lemma}\label{lem:diam}
	Given a matrix $A\in\reals^{n \times p}$, denote
	\BEQ
	\frac{1}{2}\diam(B_1^p \cap \Null(A))  = \sup_{\substack{Az=0\\\|z\|_1\leq 1}} \|z\|_2,
	\EEQ
	the radius of the section of the $\ell_1$ ball $B_1^p$ by the null space of the matrix $A$ and 
	\BEQ\label{eq:s_A}
	s_A  \triangleq 1/\diam(B_1^p \cap \Null(A))^2,
	\EEQ
	a recovery threshold.
	Then $A$ satisfies \eqref{def:nsp} at any order $s < s_A$ with constant
	\BEQ\label{eq:alpha_recovery_threshold}
	\alpha = 2\sqrt{s_A/s} -1 >1.
	\EEQ
\end{lemma}
\begin{proof}
	For any $z\in\Null(A)$ and support set $S$ with $\Card(S)\leq s$, using equivalence of norms and definition of the radius,
	\BEQ
	\|z_S\|_1 \leq \sqrt{s} ~\|z\|_2 \leq \frac{1}{2}\sqrt{\frac{s}{s_A}}~ \|z\|_1 = \frac{1}{2}\sqrt{\frac{s}{s_A}}~(\|z_S\|_1+\|z_{S^c}\|_1),
	\EEQ
	which means, as $s <s_A$,
	\BEQ
	\|z_{S^c}\|_1\geq (2\sqrt{s_A/s} -1 )\|z_S\|_1,
	\EEQ
	hence the desired result.
\end{proof}

With $s_A$ defined in \eqref{eq:s_A}, for any signal $x^*$ of sparsity $s <s_A$, the sharpness bound~\eqref{eq:sharpness} then reads 
\BEQ
\|x\|_1-\|x^*\|_1 \geq \left(1-\sqrt{s/s_A}\right) \|x-x^*\|_1,
\EEQ
and the optimal restart scheme defined in Proposition~\ref{prop:complex-l1} has complexity
\BEQ
\|\hat y\|_1-\|x^*\|_1 \leq \exp\left( -\left(1-\sqrt{s/s_A} \right)\frac{e}{ 2\sqrt{p}}  N  \right) \epsilon_0,
\EEQ
which means that, given a sensing matrix $A$ with recovery threshold $s_A$, the sparser the signal, the faster the algorithm. 

Precise estimates of the diameter of random sections of norm balls can be computed using classical results in geometric functional analysis. The low $M^*$ estimates of \cite{Pajo86} (see \cite[Theorem 3.1]{Vers11} for a concise presentation) show that when $E\subset \reals^p$ is a random subspace of codimension $n$ (e.g. the null space of a random matrix~$A\in\reals^{n\times p}$), then
\BEQ
\diam(B_1^p \cap E)\leq c \sqrt{\frac{\log p}{n}},
\EEQ
with high probability, where $c>0$ is an absolute constant. This means that the recovery threshold $s_A$ satisfies
\BEQ
s_A \geq n/(c^2\log p),
\EEQ
with high probability and leads to the following corollary.
\begin{corollary}
	Given a random sampling matrix $A\in\reals^{n\times p}$ whose nullspace is Haar distributed on the Grassman manifold, and a signal $x^*$ with sparsity 	
	$s<n/(c^2\log p)$, \eqref{eq:restart} scheme with optimal parameters defined in \eqref{eq:optimal_param} 
	outputs a point $\hat y$ such that
	\BEQ
	\|\hat y\|_1-\|x^*\|_1 \leq \exp\left( -\left(1-c\sqrt{\frac{s\log p}{n}} \right)\frac{e}{ 2\sqrt{p}}  N  \right) \epsilon_0,
	\EEQ
	with high probability, where $c$ is a universal constant and $N$ is the total number of iterations.
\end{corollary}
This means that the complexity of the optimization problem~\eqref{eq:l1-recov} is controlled by the oversampling ratio $n/s$. In other words, while increasing the number of samples increases the time complexity of elementary operations of the algorithm, it also increases its rate of convergence. 

\section{Renegar's condition number and restricted singular values}\label{s:conic_linear}
We first gave concrete evidence of the link between optimization complexity and recovery performance for the exact recovery problem by highlighting sharpness properties of the objective around the true signal, given by the null space condition. We now take a step back and consider results on the underlying conic geometry of recovery problems that also control both computational and statistical aspects. 

On the statistical side, minimal conically restricted singular values are known to control recovery performance in robust recovery problems. On the computational side, Renegar's condition number, a well known computational complexity measure for conic convex programs, controls the cost of obtaining optimality certificates for exact recovery and the sharpness of exact recovery problems (hence computational complexity of the \eqref{eq:restart} scheme presented in the previous section). 
Numerical experiments will then illustrate its relevance to control numerous other classical algorithms. By observing that minimal conically restricted singular values match the worst case of Renegar's condition number on sparse signals, our analysis shows once more that one single geometrical quantity controls both statistical robustness and computational complexity of recovery problems.

\subsection{Conic linear systems}\label{s:conic_results}
Conic linear systems arise naturally from optimality conditions of the exact recovery problem. To see this, define the tangent cone at point $x$ with respect to the $\ell_1$ norm, that is, the set of descent directions for $\|\cdot\|_1$ at $x$, as
\BEQ\label{eq:tan}
\mathcal{T}(x)= \mathrm{cone}\{z: \|x+z\|_1\leq \|x\|_1\}.
\EEQ
As shown for example by \cite[Prop\,2.1]{Chan12} a point $x$ is then the unique optimum of the exact recovery problem~\eqref{eq:l1-recov} if and only if $\Null(A)\cap \mathcal{T}(x) = \{ 0\}$, that is, there is no point satisfying the linear constraints that has lower $\ell_1$ norm than $x$.
Correct recovery of an original signal $x^*$ is therefore certified by the infeasibility of a conic linear system of the form
\BEQ
\label{eq:coneP}
\tag{P\textsubscript{A,C}}
\BA{ll}
\mbox{find} & z\\
\mbox{s.t.} & Az=0\\
& z \in C,~ z\neq 0,
\EA
\EEQ
where $C$ is a closed cone and $A$ a given matrix. For both computational and statistical aspects we will be interested in the distance to feasibility. On the computational side this will give a distance to ill-posedness that plays the role of a condition number. On the statistical side it will measure the amount of perturbation that the recovery can handle. 

\begin{definition}[Distance to feasibility] Writing $\mathcal{M}_{C} = \{ A \in \reals^{n\times p} : \eqref{eq:coneP} \textrm{ is infeasible} \}$, distance to feasibility is defined as
\BEQ\label{def:distfeas}
\sigma_C(A) \triangleq \inf_{\Delta A } ~\{ \|\Delta A \|_2 : A + \Delta A  \notin \mathcal{M}_C \}.
\EEQ
\end{definition}
A geometric analysis of the problem explicits the distance to feasibility in terms of minimal conically restricted singular value, as recalled in the following lemma.
\begin{lemma}\label{lem:dist_feas_conic_eig}
	Given a matrix $A\in \reals^{p\times n}$ and a closed cone $C$, the distance to feasibility of \eqref{eq:coneP} is given by
	\BEQ\label{def:mcrsv}
	\sigma_C(A) = \min_{\substack{x \in C \\ \|x\|_2 =1}} \|Ax\|_2.
	\EEQ
\end{lemma}
\begin{proof}
	We recall the short proof of \cite[Lemma \,3.2]{Amel14a}. Similar results have been derived by \cite[Theorem \,2]{Freu99} and \cite[Lemma \,3.2]{Bell09}.
	Let $z \in C$, with $\|z\|_2 = 1$, achieve the minimum above. 
	Then $\Delta A  = -Azz^T$ satisfies $(A+\Delta A )z =0$, so $A+\Delta A  \notin \mathcal{M}_C$ and 
	\BEQ
	\sigma_C(A) \leq \|\Delta A \|_2 = \|Az\|_2\|z\|_2 = \min_{\substack{x \in C \\ \|x\|_2 =1}} \|Ax\|_2.
	\EEQ
	On the other hand denote $\Delta A $ a perturbation such that $A+\Delta A  \notin \mathcal{M}_C$. Then there exists $z\in C\setminus{\{0\}}$ such that $(A+\Delta A )z = 0$. Thus we have
	\BEQ
	\|\Delta A \|_2 \geq \frac{\|\Delta A  z\|_2}{\|z\|_2} = \frac{\|A z\|_2}{\|z\|_2} \geq \min_{\substack{x \in C \\ \|x\|_2 =1}} \|Ax\|_2.
	\EEQ
	Taking the infimum on the left-hand side over all $\Delta A$ such that $A +\Delta A \notin \mathcal{M}_C$ concludes the proof.
\end{proof}

Expression~\eqref{def:mcrsv} writes distance to infeasibility as a cone restricted eigenvalue. Minimal cone restricted eigenvalues also directly characterize recovery performance as we recall now.

\subsection{Recovery performance of robust recovery}\label{s:recov}
Several quantities control the stability of sparse recovery in a noisy setting, with e.g. \citep{Cand06a} using restricted isometry constants, and \citep{Kash07,Judi08} using diameters with respect to various norms. In this vein, the previous section showed that recovery of a signal $x^*$ is ensured by infeasiblity of the conic linear system (\textcolor{ddarkbrown}{P\textsubscript{A,$\mathcal{T}(x^*)$}}), i.e. positiveness of the minimal conically restricted singular value $\sigma_{\mathcal{T}(x^*)}(A)$. We now show how this quantity also controls recovery performance in the presence of noise. 
In that case, the robust recovery problem attempts to retrieve an original signal $x^*$ by solving
\BEQ\label{eq:app-recov}\tag{Robust $\ell_1$ recovery}
\BA{ll}
\mbox{minimize} & \|x\|_1\\
\mbox{subject to} & \|Ax-b\|_2\leq \delta \|A\|_2,
\EA
\EEQ
in the variable $x\in\reals^p$, with the same design matrix $A\in\reals^{n \times p}$, where $b\in\reals^n$ are given observations perturbed by noise of level $\delta>0$. The following classical result then bounds reconstruction error in terms of $\sigma_{\mathcal{T}(x^*)}(A)$.
\begin{lemma}\label{lem:error-bnd}
	Given a coding matrix $A \in \reals^{n \times p}$ and an original signal $x^*$, suppose we observe $b=Ax^*+w$ where $\|w\|_2\leq \delta \|A\|_2$ and denote an optimal solution of~\eqref{eq:app-recov} by $\hat x$.
	If the minimal singular value $\sigma_{\mathcal{T}(x^*)}(A)$ in~\eqref{def:mcrsv} restricted to the tangent cone $\mathcal{T}(x^*)$ in~\eqref{eq:tan} is positive, the following error bound holds:
	\BEQ\label{eq:error-bnd}
	\|\hat x-x^*\|_2\leq 2\frac{\delta\|A\|_2}{\sigma_{\mathcal{T}(x^*)}(A)}.
	\EEQ
\end{lemma}
\begin{proof}
	We recall the short proof of \cite[Prop.\,2.2]{Chan12}. Both $\hat x$ and $x^*$ are feasible for \eqref{eq:app-recov} and $\hat x$ is optimal, so that $\|\hat x\|_1 \leq \|x^*\|_1$. Thus, the error vector $\hat x - x^*$ is in the tangent cone $\mathcal{T}(x^*)$. By the triangle inequality,
	\begin{align}
	\|A(\hat x - x^*)\|_2 & \leq \|A\hat x - b\|_2 + \|Ax^* - b\|_2 \leq 2\delta \|A\|_2.
	\end{align}
	Furthermore, by definition of $\sigma_{\mathcal{T}(x^*)}(A)$,
	\begin{align}
	\|A(\hat x - x^*)\|_2 & \geq \sigma_{\mathcal{T}(x^*)}(A)\,  \|\hat x - x^*\|_2.
	\end{align}
	Combining the two concludes the proof.
\end{proof}

Therefore the robustness of the coding matrix $A$ on all $s$-sparse signals is controlled by 
\BEQ\label{def:mu_s}
\mu_s(A) \triangleq  \inf_{x \,:\, \|x\|_0 \leq s }~ \min_{\substack{z\in\mathcal{T}(x) \\ \|z\|_2 = 1}} \|Az\|_2.
\EEQ
Expression of this minimal conically singular value can the be simplified by identifying the tangent cones on $s$-sparse signals, as done in the following lemma.
\begin{lemma}\label{lem:mu=nu}
	For any subset $S \subset \llbracket 1, p \rrbracket$, let
	\BEQ
	\mathcal{E}_S = \{z : \|z_{S^c}\|_1 \leq \|z_S\|\}\qquad \text{and} \qquad \mathcal{F}_S = \bigcup_{x \, : \,  x = x_S} \mathcal{T}(x),
	\EEQ
	then $\mathcal{E}_S = \mathcal{F}_S$.
\end{lemma}
\begin{proof}
	Let $z \in \mathcal{E}_S$, take $x = -z_S$, then
	\BEQ
	\|x+z\|_1 = \|z_{S^c}\|_1 \leq \|z_S\|_1 = \|x\|_1.
	\EEQ
	Therefore $z \in \mathcal{T}(x)\subset \mathcal{F}_S$ as $x = x_S$.
	
	Conversely let $z\in \mathcal{F}_S$, and $x \in \reals^p$, with $x =x_S$, such that $z\in \mathcal{T}(x)$. Then
	\BEQ
	\|x+z\|_1 = \|x+z_S\|_1 + \|z_{S^c}\|_1 \geq \|x\|_1 -\|z_S\|_1+\|z_{S^c}\|_1.
	\EEQ
	As $z \in \mathcal{T}(x)$, this implies $\|z_{S^c}\|_1 \leq \|z_S\|_1$, so $z\in\mathcal{E}_S$ and we conclude that $\mathcal{E}_S = \mathcal{F}_S$.
\end{proof}

Therefore, the previous expression for the minimal conically restricted singular value~\eqref{def:mu_s} can be equivalently stated as
\BEQ\label{def:mu_s_sparse}
\mu_s(A) = \min_{\substack{S\subset\llbracket 1,p \rrbracket \\ \Card(S) \leq s}} ~\min_{\substack{\|z_{S^c}\|_1\leq \|z_S\|_1 \\ \|z\|_2 = 1}} \|Az\|_2.
\EEQ
This quantity upper bounds the conically restricted singular value introduced in \citep{Bick07} defined as
\[
\kappa_s(A) = \inf_{\substack{S \in \llbracket 1, p \rrbracket \\ \Card(S) \leq s}} \inf_{\substack{\|x_{S^c}\|_1 \leq  \|x_S\|_1 \\ \|x_S\|_2 =1}}  \|Ax\|_2 
\]
In \cite{Bick07}, the authors showed that it controls estimation performance of LASSO and Dantzig selector, which was further explored by \citep{Van-09}. \. Observe that positiveness of $\mu_s(A)$ is equivalent to \eqref{def:nsp} at order $s$ with constant $1$ which shows necessity of \eqref{def:nsp} for sparse recovery. 

Since both null space property and conically restricted singular values are necessary and sufficient conditions for exact recovery they may have been linked previously in the literature. Here we derive estimates for the constant in~\eqref{def:nsp} from the minimal cone restricted singular value using tools form conic linear systems. We search for $\alpha$ such that~\eqref{def:nsp} is satisfied at order $s$. Equivalently we search for $\alpha$ such that for any support $S$ of cardinality at most $s$, the conic linear system 
\BEQ\label{eq:conic_nsp_alpha}
\BA{ll}
\mbox{find} & z\\
\mbox{s.t.} & Az=0\\
& \|z_{S^c}\|_1 \leq \alpha \|z_S\|_1,~ z\neq 0
\EA
\EEQ
is infeasible. Notice that system~\eqref{eq:conic_nsp_alpha} for $\alpha>1$ is a perturbed version of the case $\alpha=1$, so the problem reduces to studying the sensitivity to perturbations of conic linear systems as shown in the following lemma.

\begin{lemma}\label{lem:nsp_mu}
	Given a matrix $A\in \reals^{n \times p}$ and an integer $s \in \llbracket 1, p\rrbracket$,
	if the minimal conically restricted singular value $\mu_s(A)$ in \eqref{def:mu_s} and  \eqref{def:mu_s_sparse} is positive, then $A$ satisfies \eqref{def:nsp} at order $s$ for any constant 
	\BEQ
	\alpha \leq \left(1- \frac{\mu_s(A)}{\|A\|_2}\right)^{-1}.
	\EEQ
\end{lemma}
\begin{proof}
	For a support $S$ of cardinality at most $s$, write $P$ the orthogonal projector on this support (that is, $Px = x_S$), $\bar P = \idm-P$ its orthogonal projector and define the closed cone $C_S = \{  z : \|z_{S^c} \|_1 \leq \|z_S\|_1\}$. Given $\alpha \geq 1$, denote $H = \alpha^{-1}P +\bar P = \idm -(1-\alpha^{-1})P $. Observe that
	\BEQ
	\{ z : \|z_{S^c} \|_1 \leq \alpha \|z_S\|_1 \} = H C_S.
	\EEQ
	Therefore, the conic linear system~\eqref{eq:conic_nsp_alpha} reads
	\BEQ
	\BA{ll}
	\mbox{find} & z\\
	\mbox{s.t.} & Az=0\\
	& z \in H C_S,~ z\neq 0.
	\EA
	\EEQ
	As $H$ is invertible, this is equivalent to
	\BEQ\label{eq:conic_nsp_alpha_proof}
	\BA{ll}
	\mbox{find} & z\\
	\mbox{s.t.} & AHz = 0\\
	& z \in C_S,~ z\neq 0.
	\EA
	\EEQ
	Therefore, if the conic linear system
	\BEQ
	\BA{ll}
	\mbox{find} & z\\
	\mbox{s.t.} & Az = 0\\
	& z \in C_S,~ z\neq 0
	\EA
	\EEQ
	is infeasible, that is $\sigma_{C_S}(A)>0$, by Lemma~\ref{lem:dist_feas_conic_eig}, which is true for $\mu_s >0$, then by definition of the distance to feasibility, \eqref{eq:conic_nsp_alpha_proof} is also 
	infeasible provided $\|AH-A\|_2 \leq \sigma_{C_S}(A)$, which holds for any $\alpha\geq 1$ such that
	\BEQ
	(1-\alpha^{-1})\|AP\|_2 \leq \sigma_{C_S}(A).
	\EEQ
	Using that $\|AP\|_2 \leq \|A\|_2$, infeasibility is ensured in particular for any $\alpha$ such that
	\BEQ
	1-\frac{\sigma_{C_S}(A)}{\|A\|_2} \leq \alpha^{-1}.
	\EEQ
	To ensure infeasibility of the conic linear systems~\eqref{eq:conic_nsp_alpha} for any support $S$, it suffices to take $\alpha$ such that
	\BEQ
	1-\frac{\mu_s(A)}{\|A\|_2} \leq \alpha^{-1}.
	\EEQ
	This means that \eqref{def:nsp} at order $s$ is satisfied for any 
	\BEQ
	\alpha \leq \left(1- \frac{\mu_s(A)}{\|A\|_2}\right)^{-1}
	\EEQ
	where we used that, by definition of the minimal conically restricted singular value, $\mu_s(A)\leq\|A\|_2$ (in case of equality \eqref{def:nsp}, will be satisfied for any $\alpha \geq 1$).
\end{proof}

We now relate the minimal cone restricted singular value to computational complexity measures.

\subsection{Computational complexity of recovery problems}{\label{s:renegar}
Computational complexity for convex optimization problems is often described in terms of polynomial functions of problem size. This produces a clear link between problem structure and computational complexity but fails to account for the nature of the data. If we use linear systems as a basic example, unstructured linear systems of dimension $n$ can be solved with complexity $O(n^3)$ regardless of the matrix values, but iterative solvers will converge much faster on systems that are better conditioned. The seminal work of \cite{Rene95,Rene01} extends this notion of conditioning to optimization problems, producing data-driven bounds on the complexity of solving conic programs, and showing that the number of outer iterations of interior point algorithms increases as the distance to ill-posedness decreases. 

\subsubsection{Renegar's condition number}
Renegar's condition number \citep{Rene95,Rene95a,Pena00} provides a data-driven measure of the complexity of certifying infeasibility of a conic linear system of the form presented in \eqref{eq:coneP} (the larger the condition number, the harder the problem). It is rooted in the sensible idea that certifying infeasibility is easier if the problem is far from being feasible. It is defined as the scale invariant reciprocal of the distance to feasibility $\sigma_C(A)$, defined in \eqref{def:distfeas}, of problem \eqref{eq:coneP}, i.e.
\BEQ\label{def:renegar}
\mathcal{R}_{C}(A) \triangleq \frac{\|A\|_2}{\sigma_{C}(A)} = \|A\|_2/\min_{\substack{x \in C \\ \|x\|_2 =1}} \|Ax\|_2.
\EEQ
Notice that, if $C$ were the whole space $\reals^p$, and if $A^TA$ were full-rank (never the case if $n < p$), then $\sigma_{C}(A)$ would be the smallest singular value of $A$. As a result, $\mathcal{R}_{C}(A)$ would reduce to the classical condition number of $A$ (and to $\infty$ when $A^TA$ is rank-deficient). Renegar's condition number is necessarily smaller (better) than the latter, as it further incorporates the notion that $A$ need only be well conditioned along those directions that matter with respect to $C$.

\subsubsection{Complexity of certifying optimality}
In a first step, we study the complexity of the oracle certifying optimality of a candidate solution $x$ to~\eqref{eq:l1-recov} as a proxy for the problem of computing an optimal solution to this problem. 
As mentioned in Section~\ref{s:conic_results}, optimality of a point $x$ is equivalent to infeasibility of
\BEQ
\label{eq:coneT}
\tag{P\textsubscript{A,$\mathcal{T}(x)$}}
\BA{ll}
\mbox{find} & z\\
\mbox{s.t.} & Az=0\\
& z \in \mathcal{T}(x),~ z\neq 0,
\EA
\EEQ
where the tangent cone $\mathcal{T}(x)$ is defined in~\eqref{eq:tan}.
By a theorem of alternative, infeasibility of \eqref{eq:coneT} is equivalent to feasibility of the dual problem
\BEQ
\label{eq:coneT_dual}\tag{D\textsubscript{A,$\mathcal{T}(x)$}}
\BA{ll}
\mbox{find} & y\\
\mbox{s.t.} &  A^Ty \in \intr(\mathcal{T}(x)^\circ),
\EA
\EEQ
where $\mathcal{T}(x)^\circ$ is the polar cone of $\mathcal{T}(x)$. Therefore, to certify infeasibility of \eqref{eq:coneT} it is sufficient to exhibit a solution for the dual problem \eqref{eq:coneT_dual}.

Several references have connected Renegar's condition number and the complexity of solving such conic linear systems using various algorithms \citep{Rene95,Freu99a,Epel00,Rene01,Vera07,Bell09a}. In particular, \cite{Vera07} linked it to the complexity of solving the primal dual pair~\eqref{eq:coneT}--\eqref{eq:coneT_dual} using a barrier method. They show that the number of outer barrier method iterations grows as 
\BEQ
O\left(\sqrt{\rho} \log\left(\rho\,\mathcal{R}_{\mathcal{T}(x)}(A)\right)\right),
\EEQ
where $\rho$ is the barrier parameter, while the conditioning (hence the complexity) of the linear systems arising at each interior point iteration is controlled by $\mathcal{R}_{\mathcal{T}(x)}(A)^2$. This link was also tested empirically on linear programs using the NETLIB library of problems by \cite{Ordo03}, where computing times and number of iterations were regressed against estimates of the condition number computed using the approximations for Renegar's condition number detailed by \cite{Freu03a}.

Studying the complexity of computing an optimality certificate gives insights on the performance of oracle based optimization techniques such as the ellipsoid method. We now show how Renegar's condition also controls the number steps in the \eqref{eq:restart} scheme presented in Section~\ref{s:restart}.

\subsubsection{Complexity of restart scheme with Renegar's condition number}
Convergence of the \eqref{eq:restart} scheme presented in Section~\ref{s:restart} is controlled by the sharpness of the problem deduced from \eqref{def:nsp}. We now observe that sharpness is controlled by the worst case Renegar condition number for the optimality certificates \eqref{eq:coneT} on all $s$-sparse signals, defined as  
\BEQ\label{def:worse_renegar}
\mathcal{R}_s(A) \triangleq \sup_{x \, : \, \|x\|_0 \leq s} \mathcal{R}_{\mathcal{T}(x)}(A)  = \|A\|_2/\mu_s(A).
\EEQ
Connecting Lemmas~\ref{lem:error-bnd},~\ref{lem:nsp_mu} and Proposition~\ref{prop:complex-l1} we get the following corollary.
\begin{corollary}
Given a coding matrix $A\in\reals^{n \times p}$ and a sparsity level $s\geq 1$, if $\mathcal{R}_s(A)<+\infty$ in~\eqref{def:worse_renegar} then optimal \eqref{eq:restart} scheme achieves an $\epsilon$ precision in at most
\BEQ
O( (2\mathcal{R}_s(A) -1) \log \epsilon^{-1})
\EEQ
iterations.
\end{corollary}
This shows that Renegar's condition number explicitly controls the convergence of an algorithmic scheme devoted to the exact recovery problem \eqref{eq:l1-recov}, through its link with sharpness.

On the statistical side, we observed that the minimal conically restricted singular value controls recovery performance of robust procedures and that its positivity ensures exact recovery. On the computational side, we presented the role of Renegar's condition number as a computational complexity measure for sparse recovery problems. A key observation is that the worst case of Renegar's condition number $\mathcal{R}_s(A)$, defined in \eqref{def:worse_renegar}, matches the minimal conically restricted singular value defined in \eqref{def:mu_s}. Once again, a single quantity controls both aspects. This at least partially explains the common empirical observation (see, e.g., \cite{Dono06}) that problem instances where statistical estimation succeeds are computationally easier to solve.

\subsection{Computational complexity for inexact recovery}
When the primal problem~\eqref{eq:coneT} is feasible, so that $\sigma_{\mathcal{T}(x)}(A) = 0$, Renegar's condition number as defined here is infinite. While this correctly captures the fact that, in that regime, statistical recovery does not hold, it does not properly capture the fact that, when~\eqref{eq:coneT} is ``comfortably" feasible, certifying so is easy, and algorithms terminate quickly (although they return a useless estimator). From both a statistical and a computational point of view, the truly delicate cases correspond to problem instances for which both~\eqref{eq:coneT} and~\eqref{eq:coneT_dual} are only barely feasible or infeasible. This is illustrated in simple numerical example by \cite[\S11.4.3]{Boyd03} and in our numerical experiments, corresponding to the peaks in the CPU time plots of the right column in Figure~\ref{fig:XP}: problems where sparse recovery barely holds/fails are relatively harder. For simplicity, we only focused here on distance to feasibility for problem~\eqref{eq:coneT}. However, it is possible to symmetrize the condition numbers used here as described by \cite[\S1.3]{Amel14a}, where a symmetric version of the condition number is defined as
\BEQ
\bar{\mathcal{R}}_{\mathcal{T}(x)}(A)=\min\left\{\frac{\|A\|}{\sigma^P_{\mathcal{T}(x)}(A)},\frac{\|A\|}{\sigma^D_{\mathcal{T}(x)}(A)}\right\},
\EEQ
where $\sigma^P_{\mathcal{T}(x)}(A)$ and $\sigma^D_{\mathcal{T}(x)}(A)$ denote the distance to feasibility of respectively \eqref{eq:coneT} and \eqref{eq:coneT_dual}.
This quantity peaks for programs that are nearly feasible/infeasible. 

As we noticed in Section~\ref{s:restart}, a {\L}ojasiewicz inequality~\eqref{eq:Loja} for the \eqref{eq:l1-recov} problem is sufficient to ensure linear convergence of the restart scheme. Connecting the symmetrized Renegar condition number to the {\L}ojasiewicz inequality constant $\gamma$ may then produce complexity bounds for the restart scheme beyond the recovery case.
{\L}ojasievicz inequalities for convex programs have indeed proven their relevance. They were used by \cite{Ferc16,Roul17} to accelerate classical methods, in particular on the LASSO problem. Lower computational bounds for the computational complexity of accelerated methods on convex optimization problems satisfying sharpness assumptions were also studied by \cite[Page 6]{Nemi85}. Although the {\L}ojasievicz inequality is proven to be satisfied by a broad class of functions \citep{Bolt07}, quantifying its parameters is still a challenging problem that would enable better parameter choices for appropriate algorithms.


\subsection{Other algorithms}
The restart scheme presented in Section~\ref{s:restart} is of course not the only one to solve problem~\eqref{eq:l1-recov} in practice and it has not been analyzed in the noisy case.  
However, we will observe in the numerical experiments of Section~\ref{s:numres} that the condition number is correlated with the empirical performance of efficient recovery algorithms such as LARS \citep{Efro04} and Homotopy~\citep{Dono06,Asif14}. On paper, the computational complexities of~\eqref{eq:l1-recov} and~\eqref{eq:app-recov} are very similar (in fact, infeasible start primal-dual algorithms designed for solving~\eqref{eq:l1-recov} actually solve problem \eqref{eq:app-recov} with $\delta$ small). 
However in our experiments, we did observe sometimes significant differences in behavior between the noisy and noiseless case.

\section{Generalization to Common Sparsity Inducing Norms}\label{s:beyond_l1}
In this section we generalize previous results to sparse recovery problems in (non-overlapping) group norms or nuclear norm. Group norms arise in contexts such as genomics to enforce the selection of groups of genes (e.g.,~\cite{Oboz11} and references therein.) The nuclear norm is used for low-rank estimation (e.g.,~\cite{Rech08} and references therein.) 
We use the framework of decomposable norms introduced by \cite{Nega09} which applies to these norms. This allows us to generalize the null space property and to derive corresponding sharpness bounds for the exact recovery problem in a broader framework. We then again relate recovery performance and computational complexity of these recovery problems.

\subsection{Decomposable norms}\label{s:decomposable}
Sparsity inducing norms have been explored from various perspectives. Here, we use the framework of decomposable norms by \cite{Nega09} to generalize our results from $\ell_1$ norms to non-overlapping group norms and nuclear norms in a concise form. We then discuss the key geometrical properties of these norms and potential characterization of their conic nature.

We first recall the definition of decomposable norms by \cite{Nega09} in terms of projectors.
\begin{definition}{\bf Decomposable norms}\label{def:decomposable}
	Given a Euclidean space $E$, a norm $\|.\|$ on $E$ is said to be decomposable if there exists a family of orthogonal projectors $\mathcal{P}$ such that
	\begin{enumerate}[label=(\roman*)]
		\item to each $P\in \mathcal{P}$ is associated a non-negative weight $\eta(P)$ and an orthogonal projector $\bar{P}$ such that $P\bar P = \bar P P = 0$, and
		\item for any $x \in E$ and $P \in \mathcal{P}$, $\|Px + \bar P x\| = \|Px\| + \|\bar Px\|$ \label{eq:decomposability}.
	\end{enumerate}
	A signal $x$ is then said to be $s$-sparse if there exists $P \in \mathcal{P}$, such that $\eta(P)\leq s$ and $Px =x$.
\end{definition}

We now detail the family of projectors for some decomposable norms of interest.
\subsubsection{$\ell_1$ norm} 
In the the $\ell_1$ norm case, $E = \reals^p$ and $\mathcal{P}$ is the set of projectors on coordinate subspaces of $\reals^p$, that is, $\mathcal{P}$ contains all projectors which zero out all coordinates of a vector except for a subset of them, which are left unaffected. The maps $\bar P$ are the complementary projectors: $\bar{P} = \idm-P$. Property~\ref{eq:decomposability} is the classical decomposability of the $\ell_1$ norm. Naturally, the complexity level corresponds to the number of coordinates preserved by $P$, i.e., $\nu(P) = \Rank(P)$. These definitions recover the usual notion of sparsity.

\subsubsection{Group norms} 
Given a partition $G$ of $\llbracket 1,p\rrbracket$ in (non-overlapping) groups $g \subset \llbracket 1, p \rrbracket$, the group norm is defined for $x \in \reals^p$ as
\BEQ
\|x\| = \sum_{g\in G} \|x_g\|_r,
\EEQ
where $\|x_g\|_r$ is the $\ell_r$-norm of the projection of $x$ onto the coordinates defined by $g$. The cases $r =2,\infty$ correspond respectively to $\ell_1/\ell_2$ and $\ell_1/\ell_\infty$ block norms. Here, $E = \reals^p$ and the family $\mathcal{P}$ is composed of orthogonal projectors onto coordinates defined by (disjoint) unions of groups $g$, and $\bar P = \idm-P$. Formally, to each $P$ we associate $ F \subset G$ such that for any $x\in E$, $(Px)_g =x_g$ if $g \in F$ and $(Px)_g = 0$ otherwise.  Decomposability~\ref{eq:decomposability} then clearly holds.
To each group $g$ we associate a weight $\eta_g$ and for a projector $P \in \mathcal{P}$ with associated $ F \subset G$, $\eta(P) = \sum_{g \in F} \eta_g$. A classical choice of weights is $\eta_g = 1$ for all $g\in G$.

\subsubsection{Nuclear norm}
The nuclear norm is defined for matrices $X \in \reals^{p \times q}$ with singular values $\sigma_i(X)$ as 
\BEQ
 \|X\| = \sum_{k=1}^{\min(p,q)} \sigma_k(X).
\EEQ
Here $E = \reals^{p \times q}$ and its associated family of projectors contains $P$ such that 
\BEQ
P \colon X \mapsto P_\mathrm{left}XP_\mathrm{right},
\EEQ
and
\BEQ
\bar{P} \colon X \mapsto (\idm-P_\mathrm{left})X(\idm-P_\mathrm{right}),
\EEQ 
where $P_\mathrm{left} \in \reals^{p\times p}$ and $P_\mathrm{right} \in \reals^{q\times q}$ are orthogonal projectors. Their weights are defined as $\eta(P) = \max\left(\Rank(P_\mathrm{left}),\Rank(P_\mathrm{right})\right)$ defining therefore $s$-sparse matrices as matrices of rank at most $s$. As $P$ and $\bar P$ project on orthogonal row and column spaces, condition~\ref{eq:decomposability} holds.  

Decomposable norms offer a unified nomenclature for the study of sparsity inducing norms. However, they appear to be essentially restricted to the three examples presented above. Moreover, it is not clear if their definition is sufficient to characterize the conic nature of these norms, in particular in the nuclear norm case that will require additional linear algebra results. In comparison, the framework proposed by \cite{Judi14} can encompass \emph{non-latent} overlapping groups. For future use, we simplify the third property of their definition~\citep[Section~2.1]{Judi14} in Appendix~\ref{app:lemma_sparsity_struct}. It is not clear how this view can be used for latent overlapping group norms presented by \cite{Oboz11} applied in biology. Moreover the sufficient conditions that \cite{Judi14} present are sufficient but not necessary in the nuclear norm case. 
Better characterizing the key geometrical properties of these norms is therefore a challenging research direction. 

\subsection{Sharpness and generalized null space property}
From now on, we assume that we are given an ambient Euclidean space $E$ with one of the three decomposable norms $\|.\|$ presented in previous section, i.e. $\ell_1$, group or nuclear norm, and the associated family of orthogonal projectors $\mathcal{P}$ as introduced in Definition~\ref{def:decomposable}. We study the sparse recovery problem 
\BEQ\tag{Sparse recovery}\label{eq:gen_recov}
\BA{ll}
\mbox{minimize} & \|x\|\\
\mbox{subject to} & A(x)=b
\EA
\EEQ
in the variable $x\in E$, where $A$ is a linear operator onto $\reals^n$ and the observations $b \in \reals^n$ are taken from an original point $x^*$ such that $b = A(x^*)$. We begin by generalizing the null space property in this setting.
\begin{definition}{\bf (Generalized Null space Property) } A linear operator $A$ on $E$ satisfies the Generalized Null Space Property (GNSP) \emph{for orthogonal projector $P \in \mathcal{P}$} with constant $\alpha \geq 1$ if and only if for any $z\in \Null(A)\setminus\{0\}$ such that $z = Pz + \bar P z$,
	\BEQ\label{def:gen_nsp}\tag{GNSP}
	\alpha \|Pz\| < \|\bar P z\|.
	\EEQ
	The linear operator $A$ satisfies the Generalized Null Space Property \emph{at order $s$} with constant $\alpha \geq 1$ if it satisfies it for any $P$ such that $\eta(P)\leq s$.
\end{definition}
Notice that if $\bar P = \idm-P$, the condition $z = Pz + \bar P z$ is not restrictive.  However it will be useful to prove necessity of \eqref{def:gen_nsp} for the nuclear norm. In that case, observe that it is equivalent to the condition introduced by \cite{Oyma10}, i.e.
\BEQ
\forall z \in \Null(A)\setminus\{0\}, \quad \alpha \sum_{i=1}^{s} \sigma_i(z) < \sum_{i=s+1}^{\min(p,q)} \sigma_i(z),
\EEQ
where $\sigma_i(z)$ are the singular values of $z$ in decreasing order.
Notice also that we recover the classical Definition~\ref{def:nsp} in the $\ell_1$ case. The sharpness bound then easily follows if $\bar P = \idm-P$. In the case of the nuclear norm it requires additional linear algebra results.
\begin{proposition}\label{prop:gen_recov}
	Given a linear operator $A$ that satisfies \eqref{def:gen_nsp} at order $s$ with constant $\alpha$, if the original point $x^*$ is $s$-sparse, then for any $x\in E$ satisfying $A(x) = b$, $x \neq x^*$, we have 
	\BEQ
	\|x\| - \|x^*\| > \frac{\alpha-1}{\alpha+1}\|x-x^*\|.
	\EEQ 
	This implies recovery, i.e., optimality of $x^*$ for \eqref{eq:gen_recov}.
\end{proposition} 
\begin{proof}
	Denote $P$ such that $\eta(P)\leq s$ and $Px^*= x^*$, which defines its sparsity. Let $x\neq x^*$ such that $A(x) =b$, so $z = x-x^*\in \Null(A)$ and $z\neq 0$. If $\bar P = \idm-P$, $x = Px +\bar Px$ and using the decomposability~\ref{eq:decomposability}, we have
	\begin{align}
	\|x\| &= \|Px^*+Pz\| + \|\bar{P}z\| \\
	& \geq \|x^*\| -\|Pz\|+\|\bar P z\| \\
	& = \|x^*\| +\|z\| - 2\|P z\|.
	\end{align}
	By using~\eqref{def:gen_nsp}, $\|z\| = \|Pz\|+\|\bar P z\|>(1+\alpha)\| P z\|$. The result follows by arranging the terms.
	
	If $\|.\|$ is the nuclear norm and $\bar P \neq \idm-P$, as in~\citep[Lemma 6]{Oyma10}, we use that (see~\cite[Theorem 7.4.9.1]{HornJohnson90})
	\BEQ
	\|x^*+z\| \geq \sum_{i=1}^{\min(p,q)}|\sigma_i(x^*)-\sigma_i(z)|,
	\EEQ
	where $\sigma_i(x^*),\sigma_i(z)$ denote the singular values in decreasing order of respectively $x^*$ and $z$.
	Then, using that $x^*$  has rank at most $s$,
	\begin{align}
	\|x\| &\geq  \sum_{i=1}^{s}|\sigma_i(x^*)-\sigma_i(z)| + \sum_{i=s+1}^{\min(p,q)}\sigma_i(z) \\
	& \geq  \sum_{i=1}^{s}\sigma_i(x^*) - \sum_{i=1}^{s} \sigma_i(z) + \sum_{i=s+1}^{\min(p,q)}\sigma_i(z)\\
	& = \|x^*\| -\|Qz\|+\|\bar Qz\|,
	\end{align}
	where $Q$ is the projector on the $s$ largest singular directions of $z$ and therefore $\bar Q$ the projector on the $n-s$ others. These can be defined using the singular value decomposition of $z$ such that $z = Qz + \bar Q z$. Then, using \eqref{def:gen_nsp} and the decomposability~\ref{eq:decomposability} concludes the proof as above.
\end{proof}

This shows that the sharpness bound of the form~\eqref{eq:sharpness} generalizes to non-overlapping group norms and the nuclear norm. Proposition~\ref{prop:sharp_to_nsp} can also be generalized directly to this case with our definition of \eqref{def:gen_nsp}. The smoothing argument and restart schemes developed in Section~\ref{s:restart} can then be applied with similar linear convergence rates that essentially depend on the sharpness constant. By looking at the diameter of the section of the unit ball of the norm by the null space of $A$, one may also show that the oversampling ratio controls the sharpness bound as in Section~\ref{s:oversampling}.

As in Section~\ref{s:conic_linear}, we now study the conic quantities which control statistical and optimization aspects.

\subsection{Robust recovery performance and computational complexity}\label{s:gen_renegar}
In this section, for a Euclidean space $E$ and $x\in E$ we denote $\|x\|_2$  the $\ell_2$ norm of its coefficients, if $E$ is a matrix space $\|x\|_2$ is then the Frobenius norm of $x$.
\subsubsection{Generalized cone restricted singular value}
We begin by addressing the recovery performance of robust sparse recovery problems that reads 
\BEQ\label{eq:gen_app-recov}\tag{Robust sparse recovery}
\BA{ll}
\mbox{minimize} & \|x\|\\
\mbox{subject to} & \|A(x)-b\|_2\leq \delta \|A\|_2,
\EA
\EEQ
in the variable $x\in E$, with the same linear operator $A$, where the observations $b\in \reals^n$ are affected by noise of level $\delta>0$. For a linear operator $A$ from $E$ to $\reals^n$, we denote its operator norm with respect to $\|\cdot\|$, $\|A\|_2 = \sup_{x \in E : \|x\|_2\leq 1} \|A(x)\|_2$. 

The results of Section~\ref{s:recov} transpose directly to the general case by replacing $\|\cdot\|_1$ by $\|\cdot\|$. Precisely, assuming that $b=Ax^*+w$ where $\|w\|_2\leq \delta \|A\|_2$, an optimal solution $\hat{x}$ of problem~\eqref{eq:gen_app-recov} satisfies the error bound
\BEQ 
	\|\hat x-x^*\|_2\leq 2\frac{\delta\|A\|_2}{\sigma_{\mathcal{T}(x^*)}(A)},
\EEQ
where the tangent cone is defined as 
\BEQ
\mathcal{T}(x)= \mathrm{cone}\{z: \|x+z\|\leq \|x\|\},
\EEQ
and robust recovery of $s$-sparse signals is therefore controlled by
\BEQ\label{def:gen_mu_s}
\mu_s(A) = \inf_{\substack{P \in \mathcal{P} \,:\, \eta(P)\leq s}} \: \inf_{\substack{x \in E \,:\, Px = x}} \: \min_{\substack{z\in\mathcal{T}(x) \\ \|z\|_2 = 1}} \|Az\|_2.
\EEQ

The key point is then to characterize the tangent cones of $s$-sparse signals. First, this will allow statistical estimations of $\mu_s(A)$. Second, it will enable us to estimate the constant~\eqref{def:gen_nsp}, hence sharpness of the exact recovery problem and computational complexity of associated restart schemes. This is the aim of the following lemma.
\begin{lemma}\label{lem:gen_key_lemma}
	For a given sparsity $s$, write
	\BEQ
	\mathcal{E} = \bigcup_{P\in \mathcal{P} \,:\, \eta(P) \leq s} \{z \in E \,:\, z = Pz + \bar P z, \:\|\bar P z\| \leq \|Pz\|\} 
	\EEQ
	and
\BEQ	
	 \mathcal{F} = \bigcup_{P\in \mathcal{P} \,:\, \eta(P) \leq s}\: \bigcup_{x \in E \,:\,  x = Px} \mathcal{T}(x).
\EEQ	 
	Then $\mathcal{E} = \mathcal{F}$.
\end{lemma}
\begin{proof}
	Let $z \in \mathcal{E}$ and $P \in \mathcal{P}$ such that $z = Pz +\bar P z$. Taking $x = -Pz$ we get
	\BEQ
	\|x +z\| = \|\bar P z\| \leq  \|Pz\| = \|x\|.
	\EEQ
	Therefore $z \in \mathcal{T}(x)\subset \mathcal{F}$. Conversely, if $z \in \mathcal{F}$, denote $x \in E$ and $P \in \mathcal{P}$ such that $x = Px$, $z \in \mathcal{T}(x)$ and $\eta(P)\leq s$. If $\bar P = \idm-P$, by decomposability~\ref{eq:decomposability},
	\BEQ
	\|x+z\| = \|Px+Pz\| +\|\bar Pz\| \geq \|x\| - \|P z\| + \|\bar P z\|.
	\EEQ
	Since $z\in \mathcal{T}(x)$, we have $\|x+z\| \leq \|x\|$; combined with the previous statement, this implies that $z \in \{z \in E : z = Pz + \bar P z, \:\|\bar P z\| \leq \|Pz\|\} \subset \mathcal{E}$.  Now, if $\|.\|$ is the nuclear norm, as in the proof of Proposition~\eqref{prop:gen_recov}, we have
	\BEQ
	\|x+z\| \geq \|x\| -\|Qz\|+\|\bar Qz\|,
	\EEQ
	where $Q$ is the projector on the $s$ largest singular directions of $z$ given by the singular value decomposition of $z$, so that $z = Qz+\bar Q z$. Therefore, $z\in \mathcal{T}(x)$ implies $z \in \{z\in E : z = Qz+\bar Q z, \: \|\bar Q z\| \leq \|Q z\|\}\subset \mathcal{E}$. In all cases we have therefore proven $\mathcal{E}= \mathcal{F}$.
\end{proof}

Using the previous lemma, the minimal cone restricted singular value reads:
\BEQ\label{def:gen_mu_s_sparse}
\mu_s(A) = \inf_{P \in \mathcal{P}, \: \eta(P)\leq s} \min_{\substack{z \in E, \: \|z\|_2 =1\\ z = Pz +\bar Pz, \: \|\bar P z\| \leq \|Pz\|}} \|Az\|_2.
\EEQ
This quantity can then be linked to the~\eqref{def:gen_nsp} constant, as shown in the following lemma.
\begin{lemma}\label{lem:gen_nsp_mu}
Given a linear operator $A$ on $E$, 
If the minimal cone restricted singular value $\mu_s(A)$, defined in \eqref{def:gen_mu_s} and reformulated in \eqref{def:gen_mu_s_sparse}, is positive, then $A$ satisfies \eqref{def:gen_nsp} at order $s$ for any constant 
\BEQ
\alpha \leq \left(1- \frac{\mu_s(A)}{\|A\|_2}\right)^{-1}.
\EEQ
\end{lemma}
\begin{proof}
	For a given $P\in \mathcal{P}$, denote $C_P = \{z \in \Im(P)+\Im(\bar P) \,:\, \|\bar P z\| \leq \|Pz\|\}$ and define for $\alpha \geq 1$ the conic linear system
	\BEQ\label{eq:gen_conic_nsp_alpha}
	\BA{ll}
	\mbox{find} & z \in \Im(P)+\Im(\bar P)\\
	\mbox{s.t.} & A(z)=0\\
	& \|\bar Pz\| \leq \alpha\|Pz\|,~ z\neq 0.
	\EA
	\EEQ
	Infeasibility of this system for all $P\in \mathcal{P}$ such that $\eta(P) \leq s$ is then equivalent to \eqref{def:gen_nsp} at order $s$ with constant $\alpha$. Denote $H = \idm -(1-\alpha^{-1})P$ such that 
	\BEQ
	\{ z \in \Im(P)+\Im(\bar P) : \|\bar Pz\| \leq \alpha\|Pz\|\} = H C_P.
	\EEQ
	Since $H$ is invertible, we observe as in Lemma~\ref{lem:nsp_mu} that the conic linear system \eqref{eq:gen_conic_nsp_alpha} is equivalent to
	\BEQ\label{eq:gen_conic_nsp_perturbed}
	\BA{ll}
	\mbox{find} & z \in \Im(P)+\Im(\bar P)\\
	\mbox{s.t.} & A - (1-\alpha^{-1}) APz=0\\
	& z \in C_P,~ z\neq 0.
	\EA
	\EEQ
	If this problem is infeasible for $\alpha=1$, i.e., its distance to feasibility $\mu_{C_P}(A)$ defined in \eqref{def:distfeas} is positive, then \eqref{eq:gen_conic_nsp_perturbed} is infeasible for any $\alpha\geq 1$ such that
	\BEQ
	(1-\alpha^{-1}) \|AP\| \leq \mu_{C_P}(A).
	\EEQ
	Now, if $\mu_s (A) >0$
	the conic linear system \eqref{eq:gen_conic_nsp_alpha} will still be infeasible  for any 
	\BEQ
	\alpha \leq \left(1- \frac{\mu_s(A)}{\|A\|_2}\right)^{-1}.
	\EEQ
	Thus, $A$ satisfies \eqref{def:gen_nsp} at order $s$ with $\alpha$ as above.
\end{proof}

\subsubsection{Renegar's condition number}
On the computational side, denote $\mathcal{R}_{\mathcal{T}(x)}(A)$ the  Renegar condition number of the conic linear system
\BEQ
\BA{ll}
\mbox{find} & z\\
\mbox{s.t.} & A(z)=0\\
& z \in \mathcal{T}(x),~ z\neq 0,
\EA
\EEQ
and the worst-case Renegar condition number on $s$-sparse signals
\BEQ
\mathcal{R}_s(A) \triangleq \sup_{\substack{P \in \mathcal{P} \,:\, \eta(P)\leq s}} \: \sup_{\substack{x \in E \,:\, Px = x}} \mathcal{R}_{\mathcal{T}(x)}(A)  = \|A\|_2/\mu_s(A).
\EEQ
First, Renegar's condition number plays the same role as before in computing optimality certificates for the exact recovery problems. Then, combining Lemma~\ref{lem:gen_nsp_mu} and Proposition~\ref{prop:gen_recov} shows that the sharpness bound for exact recovery reads
\BEQ
\|x\| - \|x^*\| > \frac{1}{2\mathcal{R}_s(A) -1}\|x-x^*\|.
\EEQ 
This sharpness will then control linearly convergent restart schemes for the exact recovery problem.

Overall then, as established earlier in this paper, a single geometric quantity---namely, the minimal cone restricted singular value---appears to control both computational and statistical aspects. We now illustrate this statement on numerical experiments.

\section{Numerical Results} \label{s:numres}
In this section, we first test the empirical performance of restart schemes and its link with recovery performance. We then perform similar experiments on Renegar's condition number.

\subsection{Sharpness \& restart for exact recovery} 
We test the \eqref{eq:restart} scheme on $\ell_1$-recovery problems with random design matrices. Throughout the experiments, we use the NESTA code described in \citep{Beck11} as the subroutine in the restart strategy. We generate a random design matrix $A\in\reals^{n \times p}$ with i.i.d.\ Gaussian coefficients. We then normalize $A$ so that $AA^T=\idm$ (to fit NESTA's format) and generate observations $b=Ax^*$ where $x^*\in\reals^p$ is an $s$-sparse vector whose nonzero coefficients are all ones. We denote $\hat x$ the solution given by a common solver run at machine precision and plot convergence $f(x_t) -f^* = \|x_t\|_1-\|\hat x\|_1$ (scaled such that $f(x_0)-f(\hat x) =1$).

\subsubsection{Restart scheme performance}
First we compare in Figure~\ref{fig:compare_plain} the practical scheme presented in Section~\ref{s:pract_rest} with a plain implementation of NESTA without restart or continuation steps. Dimensions of the problem are $p =300$, $n = 200$ and $s= 10$.
Starting from $x_0 = A^Tb$, we use $\epsilon_0=\|x_0\|_1$ as a first initial guess on the gap and perform a grid search of step size $h=4$ for a budget of $N= 500$ iterations. The first and last schemes of the grid search were not run as they are unlikely to produce a nearly optimal restart scheme.
The grid search can be parallelized and the best scheme found is plotted with a solid red line. The dashed red line represents the convergence rate accounting for the cost of the grid search. For the plain implementation of NESTA, we used different target precisions. 
These control indeed the smoothness of the surrogate function $f_\epsilon$ which itself controls the step size of Nesterov's algorithm. Therefore a high precision slows down the algorithm. However for low precision NESTA can be faster but will not approximate well the original signal.
Also, the theoretical bound \eqref{eq:conv_rate} might be very pessimistic, as the surrogate function $f_\epsilon$ may approximate the $\ell_1$ norm for the points of interest at a much better accuracy than $\epsilon$.

\begin{figure}[ht]
	\begin{center}
		\psfrag{fmu}[b][t]{$f(x_t) -f^*$}
		\psfrag{k}[t][b]{Inner iterations}
		\includegraphics[width=0.5\textwidth]{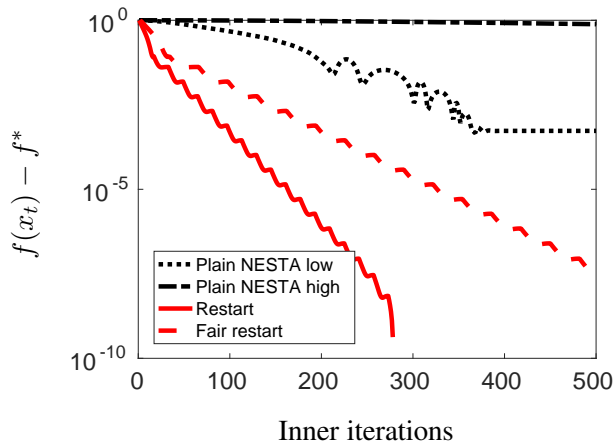} 
	\end{center}
	\caption{Best restarted NESTA (solid red line) and overall cost of the practical restart schemes (dashed red line) versus plain NESTA implementation with low accuracy $\epsilon = 10^{-1}$ (dotted black line) and higher accuracy $\epsilon = 10^{-3}$ (dash-dotted black line) for a budget of 500 iterations.}
	\label{fig:compare_plain}
\end{figure}
	
Overall, we observe a clear linear convergence of the restart scheme that outperforms the plain implementation. This was already observed by \cite{Beck11} who developed their continuation steps against which we compare in Figure~\ref{fig:compare_continuation}. We used default options for NESTA, namely 5 continuation steps with a stopping criterion based on the relative objective change in the surrogate function (specifically, the algorithm stops when these changes are lower than the target accuracy, set to $10^{-6}$). We compare continuations steps and best restart found by grid search for different dimensions of the problem, we fix $p = 300$, $s= 10$ and vary the number of samples $n= \{120,200\}$. Continuation steps converge faster with better conditioned problems, i.e., more samples. Overall the heuristic of continuation steps offer similar or better linear convergence than the restart scheme found by grid-search. Notice that a lot of parameters are involved for both algorithms, in particular the target precision may play an important role, so that more extensive experiments may be needed to refine these statements.

Our goal here is to provide a simple but strong baseline with theoretical guarantees for recovery. Improving on it, as \cite{Ferc16} did for LASSO, is an appealing research direction. Sharpness may be used for example to refine the heuristic strategy of the continuations steps.
\begin{figure}[H]
	\begin{center}
		\begin{tabular}{cc}
		\psfrag{fmu}[b][t]{$f(x_t) -f^*$}
		\psfrag{k}[t][b]{Inner iterations}
		\includegraphics[width=0.45\textwidth]{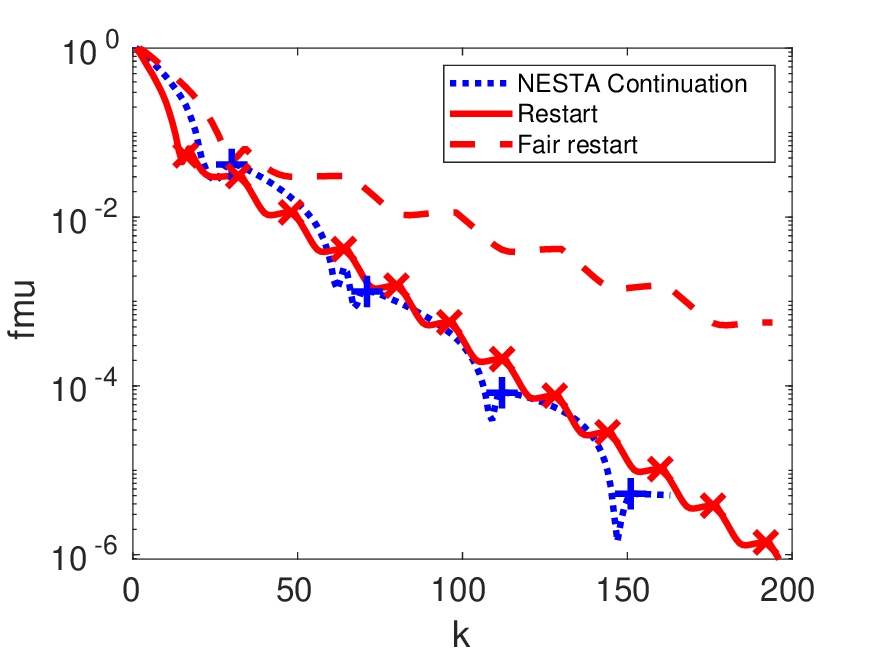} 
		&
		\psfrag{fmu}[b][t]{$f(x_t) -f^*$}
		\psfrag{k}[t][b]{Inner iterations}
		\includegraphics[width=0.45\textwidth]{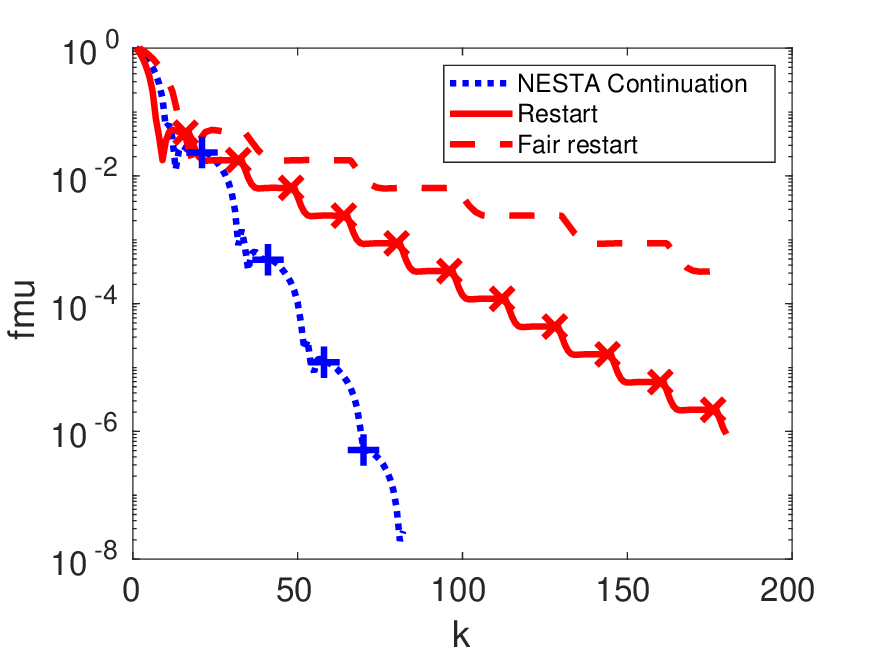} 
		\end{tabular}
	\end{center}
	\caption{Best restarted NESTA (solid red line) and overall cost of the practical restart schemes (dashed red line) versus NESTA with 5 continuation steps (dotted blue line) for a budget of 500 iterations. Crosses represent the restart occurrences. Left: $n= 120$. Right: $n = 200$.}
	\label{fig:compare_continuation}
\end{figure}

\subsubsection{Convergence rate and oversampling ratio}
We now illustrate the theoretical results of Section~\ref{s:oversampling} by running the practical scheme presented in Section~\ref{s:pract_rest} for increasing values of the oversampling ratio $\tau = n/s$. In Figure~\ref{fig:oversampling}, we plot the best scheme found by the grid search, that approximates the optimal scheme, for a budget of $N = 500$ iterations. We use a non-logarithmic grid to find the best restart scheme. Other algorithmic parameters remain unchanged: $x_0 = A^Tb$ and $\epsilon_0 = \|x_0\|_1$. We fix the dimension $p = 1000$ and either make $n$ vary for a fixed sparsity $s = 17$ or make $s$ vary for a fixed number of samples $n = 200$. These values ensure that we stay in the recovery regime as analyzed in \citep{Judi08}. In both cases we do observe an improved convergence for increasing oversampling ratio $\tau$. 

\begin{figure}[H]
\begin{center}
    \begin{tabular}{cc}
    \psfrag{fmu}[b][t]{$f(x_k)-f(x^*)$}
    \psfrag{k}[t][b]{Inner iterations}
    \psfrag{tau}{$\tau$}
    \includegraphics[width=0.45\textwidth]{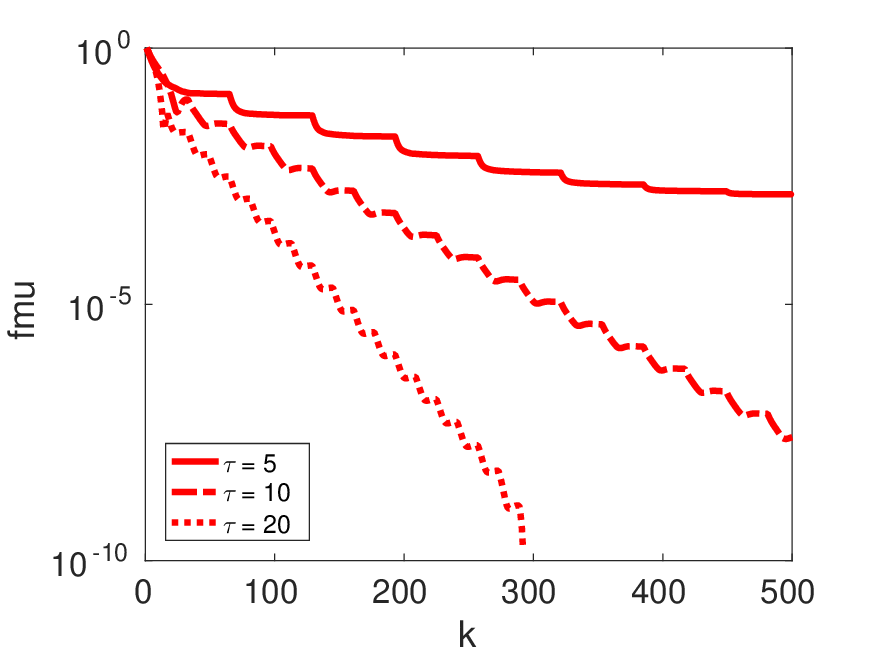}
    &
    \psfrag{fmu}[b][t]{$f(x_k)-f(x^*)$}
    \psfrag{k}[t][b]{Inner iterations}
    \psfrag{tau}{$\tau$}
    \includegraphics[width=0.45\textwidth]{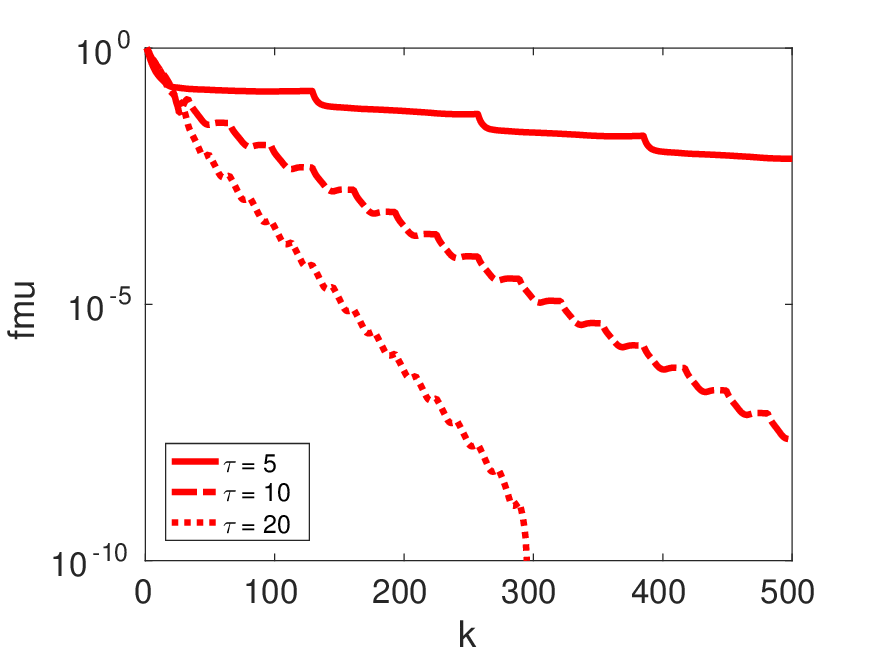}
    \end{tabular}
\end{center}
\caption{Best restart scheme found by grid search for increasing values of the oversampling ratio $\tau = n/s$ with $p=1000$. Left: sparsity $s =17$ fixed. Right: number of samples $n = 200$ fixed.}
\label{fig:oversampling}
\end{figure}



\subsection{Renegar's condition number and compressed sensing performance}
Our theoretical results showed that Renegar's condition number measures the complexity for the exact recovery problem~\eqref{eq:recov}. However it does not a priori control convergence of the robust recovery problems defined in the introduction. This numerical section aims therefore at analyzing the relevance of this condition number for general recovery problems in the $\ell_1$ case, assuming that their complexity corresponds roughly to that of checking optimality of a given point at each iteration, as mentioned in Section~\ref{s:renegar}.
We first describe how we approximate the value of $\mathcal{R}_{\mathcal{T}(x^*)}(A)$ as defined in \eqref{def:renegar} for a given original signal $x^*$ and matrix $A\in \reals^{n \times p}$. We then detail numerical experiments on synthetic data sets.

\subsubsection{Computing $\mathcal{R}_{\mathcal{T}(x^*)}(A)$}\label{s:computeC}
The condition number $\mathcal{R}_{\mathcal{T}(x^*)}(A)$ appears here in upper bounds on computational complexities and statistical performances. In order to test numerically whether this quantity truly explains those features (as opposed to merely appearing in a wildly pessimistic bound), we explicitly compute it in numerical experiments.

To compute $\mathcal{R}_{\mathcal{T}(x^*)}(A)$, we propose a heuristic which computes $\sigma_{\mathcal{T}(x^*)}(A)$ in~\eqref{def:distfeas} and~\eqref{def:mcrsv}, the value of a nonconvex minimization problem over the cone of descent directions $\mathcal{T}(x^*)$. The closure of the latter is the polar of the cone generated by the subdifferential to the $\ell_1$-norm ball at $x^*$~\citep[\S 2.3]{Chan12}. Let $S \subset \llbracket 1, p \rrbracket$ denote the support of $x^*$ and $s = \Card(S)$. Then, with $u = \operatorname{sign}(x^*)$,
\begin{align*}
\mathcal{T}(x^*) & = \operatorname{cone}\left\{ z \in \reals^p : z_S = u_S, z_{S^c} \in [-1, 1]^{p-s} ) \right\}^\circ 
= \left\{ z \in \reals^p : \|z_{S^c}\|_1 \leq - u_S^Tz_S^{} = -u^Tz \right\}.
\end{align*}
Thus, $\sigma_{\mathcal{T}(x^*)}(A)$ is the square root of
\begin{align}
\min_{z \in \reals^p} z^T A^T A z \quad \textrm{ s.t. } \quad \|z\|_2 = 1 \quad \textrm{ and } \quad \|z_{S^c}\|_1 \leq -u^Tz.
\label{eq:coneconstrainedmaxeig}
\end{align}
Let $\lambda$ denote the largest eigenvalue of $A^T A$. If it were not for the cone constraint, solutions of this problem would be the dominant eigenvectors of $\lambda \idm- A^T A$, which suggests a \emph{projected power method}~\citep{Desh14} as follows. Given an initial guess $z_{0} \in \reals^p$, $\|z_{0}\|_2 = 1$, iterate
\begin{align}
\hat z_{k+1} & = \operatorname{Proj}_{\mathcal{T}(x^*)}\left((\lambda \idm- A^T A) z_{k}\right), & z_{k+1} & = \hat z_{k+1} / \|\hat z_{k+1}\|_2,
\end{align}
where we used the orthogonal projector to $\mathcal{T}(x^*)$,
\begin{align}
\operatorname{Proj}_{\mathcal{T}(x^*)}(\tilde z) & = \arg\min_{z \in \reals^p}\|z - \tilde z\|_2^2 \quad \textrm{ s.t. } \quad \|z_{S^c}\|_1 \leq -u^Tz.
\end{align}
This convex, linearly constrained quadratic program is easily solved with CVX~\citep{Gran01}. As can be seen from KKT conditions, this iteration is a generalized power iteration~\citep{Luss13,Jour08}
\begin{align}
z_{k+1} \in \arg\max_{z\in\reals^p} z^T (\lambda \idm- A^T A) z_{k} \quad \textrm{ s.t. } \quad \|z\|_2 \leq 1 \quad \textrm{ and } \quad \|z_{S^c}\|_1 \leq -u^Tz.
\end{align}
From the latter, it follows that $\|Az_{k}\|_2$ decreases monotonically with $k$. Indeed, owing to convexity of $f(z) = \frac{1}{2}z^T(\lambda \idm- A^T A)z$, we have $f(z) - f(z_{k})\geq (z - z_{k})^T(\lambda \idm- A^T A)z_{k}$. The next iterate $z = z_{k+1}$ maximizes this lower bound on the improvement. Since $z = z_{k}$ is admissible, the improvement is nonnegative and $f(z_{k})$ increases monotonically.

Thus, the sequence $\|Az_{k}\|_2$ converges, but it may do so slowly, and the value it converges to may depend on the initial iterate $z_{0}$. On both accounts, it helps greatly to choose $z_{0}$ well. To obtain one, we modify~\eqref{eq:coneconstrainedmaxeig} by smoothly penalizing the inequality constraint in the cost function, which results in a smooth optimization problem on the $\ell_2$ sphere. Specifically, for small $\varepsilon_1, \varepsilon_2 > 0$, we use smooth proxies $h(x) = \sqrt{x^2+\varepsilon_1^2} - \varepsilon_1 \approx |x|$ and $q(x) = \varepsilon_2\log(1 + \exp(x/\varepsilon_2)) \approx \max(0, x)$. Then, with $\gamma > 0$ as Lagrange multiplier, we consider
\begin{align}
\min_{\|z\|_2 = 1} \|Az\|_2^2 + \gamma \cdot q\left(  u^T z + \sum\nolimits_{i \in S^c} h(z_i)  \right).
\end{align}
We solve the latter locally with Manopt~\citep{Boum14}, itself with a uniformly random initial guess on the sphere, to obtain $z_{0}$. Then, we iterate the projected power method. 
The value $\|Az\|_2$ is an upper bound on $\sigma_{\mathcal{T}(x^*)}(A)$, so that we obtain a lower bound on $\mathcal{R}_{\mathcal{T}(x^*)}(A)$. Empirically, this procedure, which is random only through the initial guess on the sphere, consistently returns the same value, up to five digits of accuracy, which suggests the proposed heuristic computes a good approximation of the condition number. Similarly positive results have been reported on other cones by~\cite{Desh14}, where the special structure of the cone even made it possible to certify that this procedure indeed attains a global optimum in proposed experiments. Similarly, a generalized power method was recently shown to converge to global optimizers for the phase synchronization problem (in a certain noise regime)~\citep{Boum16,zhong2017nearoptimal}. This gives us confidence in the estimates produced here.

\subsubsection{Sparse recovery performance}
We conduct numerical experiments in the $\ell_1$ case to illustrate the connection between the condition number $\mathcal{R}_{\mathcal{T}(x^*)}(A)$, the computational complexity of solving~\eqref{eq:l1-recov}, and the statistical efficiency of the estimator~\eqref{eq:app-recov}. Importantly, throughout the experiments, the classical condition number of $A$ will remain essentially constant, so that the main variations cannot be attributed to the latter.

We follow a standard setup, similar to some of the experiments by~\cite{Dono06}. Fixing the ambient dimension $p = 300$ and sparsity $s = \|x^*\|_0 = 15$, we let the number of linear measurements $n$ vary from 1 to 150. For each value of $n$, we generate a random signal $x^*\in\reals^p$ (uniformly random support, i.i.d.\ Gaussian entries, unit $\ell_2$-norm) and a random sensing matrix $A\in\reals^{n\times p}$ with i.i.d.\ standard Gaussian entries. Furthermore, for a fixed value $\delta = 10^{-2}$, we generate a random noise vector $w\in\reals^n$ with i.i.d.\ standard Gaussian entries, normalized such that $\|w\|_2 = \delta \|A\|_2$, and we let $b = Ax^* + w$. This is repeated 100 times for each value of $n$.

For each triplet $(A, x^*, b)$, we first solve the noisy problem~\eqref{eq:app-recov} with the L1-Homotopy algorithm ($\tau = 10^{-7}$)~\citep{Asif14}, and report the estimation error $\|\hat x - x^*\|_2$. Then, we solve the noiseless problem~\eqref{eq:recov} with L1-Homotopy and the TFOCS routine for basis pursuit ($\mu = 1$)~\citep{Beck12}. Exact recovery is declared when the error is less than $10^{-5}$, and we report the empirical probability of exact recovery, together with the number of iterations required by each of the solvers. The number of iterations of LARS~\citep{Efro04} is also reported, for comparison. For L1-Homotopy, we report the computation time, normalized by the computation time required for one least-squares solve in $A$, as in~\citep[Fig.\,3]{Dono06}, which accounts for the growth in $n$. Finally, we compute the classical condition number of $A$, $\kappa(A)$, as well as (a lower bound on) the cone-restricted condition number $\mathcal{R}_{\mathcal{T}(x^*)}(A)$, as per the previous section. As it is the computational bottleneck of the experiment, it is only computed for 20 of the 100 repetitions.


The results of Figure~\ref{fig:XP} show that the cone-restricted condition number explains both the computational complexity of~\eqref{eq:l1-recov} and the statistical complexity of~\eqref{eq:app-recov}: fewer samples mean bad conditioning which in turn implies high computational complexity. We caution that our estimate of $\mathcal{R}_{\mathcal{T}(x^*)}(A)$ is only a lower bound. 
Indeed, for small $n$, the third plot on the left shows that, even in the absence of noise, recovery of $x^*$ is not achieved by~\eqref{eq:app-recov}. Lemma~\ref{lem:error-bnd} then requires $\mathcal{R}_{\mathcal{T}(x^*)}(A)$ to be infinite. But the computational complexity of solving~\eqref{eq:l1-recov} is visibly favorable for small $n$, where far from the phase transition, problem~\eqref{eq:coneT} is far from infeasibility, which is just as easy to verify as it is to certify that~\eqref{eq:coneT} is infeasible when $n$ is comfortably larger than needed. This phenomenon is best explained using a symmetric version of the condition number~\citep{Amel14a} (omitted here to simplify computations).

We also solved problem~\eqref{eq:l1-recov} with interior point methods (IPM) via CVX. The number of iterations appeared mostly constant throughout the experiments, suggesting that the practical implementation of such solvers renders their complexity mostly data agnostic in the present setting. Likewise, the computation time required by L1-Homotopy on the noisy problem~\eqref{eq:app-recov}, normalized by the time of a least-squares solve, is mostly constant (at about 150). This hints that the link between computational complexity of~\eqref{eq:l1-recov} and~\eqref{eq:app-recov} remains to be fully explained.

\begin{figure}[p]
	\begin{center}
		\includegraphics[width=1.0\textwidth]{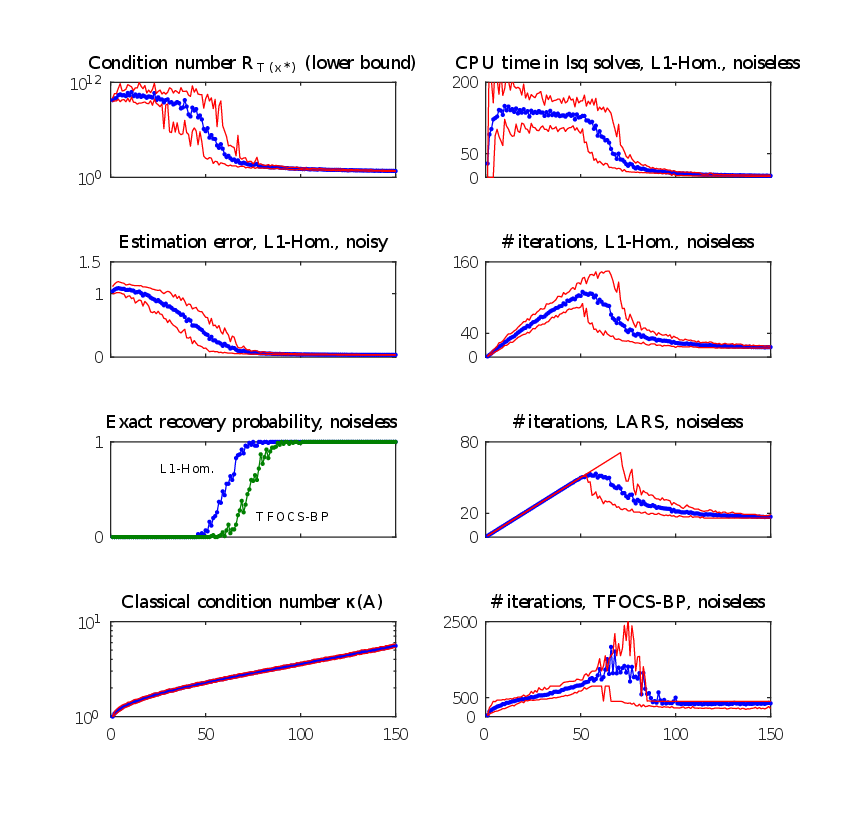}
		\caption{We plot the cone-restricted condition number of $A$ (upper left), explaining both the computational complexity of problem~\eqref{eq:l1-recov} (right column) and the statistical complexity of problem~\eqref{eq:app-recov} (second on the left). Central curves represent the mean (geometric mean in log-scale plots), red curves correspond to 10th and 90th percentile. We observe that high computing times (peaks in the right column) are directly aligned with instances where sparse recovery barely holds/fails (left), i.e. near the phase transition around $n=70$, where the distance to feasibility for problem~\eqref{eq:coneT} also follows a phase transition.}\label{fig:XP}
	\end{center}
\end{figure}

\subsection*{Acknowledgements}
VR and AA are at the D\'epartement d'Informatique at \'Ecole Normale Sup\'erieure, 2~rue Simone Iff, 75012 Paris, France. INRIA, Sierra project-team, PSL Research University. 
The authors would like to acknowledge support from a starting grant from the European Research Council (ERC project SIPA), an AMX fellowship, an NSF grant DMS-1719558, support from the {\em fonds AXA pour la recherche} and a Google focused award.

{\small{\bibliographystyle{agsm}
\bibliography{LinearCS.bib,MainPerso.bib}
}

\appendix
\section{Practical optimal restart scheme}\label{app:practical_optimal}
In Section~\ref{s:restart} we quickly give optimal restart schemes in terms of a potentially non-integer clock. Following corollary details optimal scheme for an integer optimal clock.
\begin{corollary}
		Given a coding matrix $A\in\reals^{n \times p}$ and an original signal $x^* \in \reals^p$ such that sharpness bound~\eqref{eq:sharpness} is satisfied with $\gamma >0$, running Algorithm~\ref{eq:restart} with $\rho^*$ and $t = \lceil t^* \rceil $ where $\rho^*$ and $t^*$ are defined in \eqref{eq:optimal_param}
		ensures that after $K\geq 1$ restarts, i.e. $N = K \lceil t^* \rceil$ total number of iterations,
		\BEQ
		\|\hat y\|_1-\|x^*\|_1 \leq \exp\left( -\frac{N \gamma }{2 e\sqrt{p} +\gamma}\right)  \epsilon_0.
		\EEQ
\end{corollary}
\begin{proof}
	Denote $\delta = \lceil t^*\rceil -t^* \in [0,1[$. As $\lceil t^*\rceil \geq t^*$ \eqref{eq:goal_restart} is ensured for $\rho^*$. At the $K$\textsuperscript{th} restart, $N = K(t^* + \delta)$, and
	\BEQ
	\|\hat y \|_1 -\|x^*\|_1 \leq e^{-K}\epsilon_0 = \exp( -N/(t^* + \delta)) \leq \exp(- N/(t^* + 1)).
	\EEQ
	Replacing $t^*$ by its value gives the result.
\end{proof}

\section{Remark on sparsity inducing norms}\label{app:lemma_sparsity_struct}
We quickly discuss the framework of \cite{Judi14} for sparsity inducing norms and show that it can be simplified.
We first recall the definition.
\begin{definition}\label{def:sparse_struct}{\bf (Sparsity structure \citep{Judi14}) }
	A sparsity structure on a Euclidean space $\calE$ is defined as a norm $\|\cdot\|$ on $\calE$, together with a family $\mathcal{P}$ of linear maps of $\calE$ into itself, satisfying three assumptions:
	
	\begin{enumerate}
		\item Every $P\in \mathcal{P}$ is a projector, $ P^2 = P$,
		\item Every $P \in \mathcal{P}$ is assigned a weight $\nu(P) \geq 0$ and a linear map $\bar{P}$ on $\calE$ such that $P \bar{P} =0$,
		\item \label{eq:weird_condition} For any $P\in \mathcal{P}$ and $u,v \in \calE$, one has
		\[ \|P^* u + \bar{P}^* v\|_* \leq \max (\|u\|_*,\|v\|_*), \]
	\end{enumerate}
	where $\|\cdot\|_*$ is the dual norm of $\|\cdot\|$ and $P^*$ is the conjugate mapping of the linear map $P$. 
\end{definition}

The last condition in Definition~\ref{def:sparse_struct} is arguably the least intuitive and following Lemma connects it with the more intuitive notion of decomposable norm.
\begin{lemma}\label{lem:simplification}
	Condition \eqref{eq:weird_condition} above, which reads
	\BEQ \|P^* u + \bar{P}^* v\|_* \leq \max (\|u\|_*,\|v\|_*), \EEQ
	for any $u,v \in \calE$, is equivalent to 
	\BEQ
	\|w\| \geq \|Pw\| + \|\bar P w\|,
	\EEQ
	for any $w\in \calE$.
\end{lemma}
\begin{proof}
	Denote $f : (u,v) \rightarrow \|P^* u + \bar{P}^* v\|_*$ and $g : (u,v) \rightarrow \max (\|u\|_*,\|v\|_*)$. Since $f$ and $g$ are non-negative, continuous convex functions, $f^2/2$ and $g^2/2$ are also convex continuous and following equivalences hold
	\BEQ
	f \leq g \quad \Leftrightarrow \quad \frac{f^2}{2} \leq \frac{g^2}{2} \quad \Leftrightarrow \quad  \left(\frac{f^2}{2}\right)^* \geq \left(\frac{g^2}{2}\right)^*,
	\EEQ
	using that for a convex continuous function $h$, $h^{**}= h$.
	Now combining the conjugacy result for squared norm \citep[Example 3.27]{Boyd03} showing that the conjugate of a squared norm $\|x\|^2/2$ is the squared conjugate norm $\|x\|_*^2/2$, with the result in \citep[Th.\,16.3]{Rock70}, we get
	\BEQ
	\left(\frac{f^2}{2}\right)^*(s,t) = \inf_w\{ \|w\|^2/2 : Pw=s, \bar Pw=t\},
	\EEQ
	where the infimum is $+\infty$ if the constraints are infeasible.
	Then the dual of the norm $g$ is $(s,t) \rightarrow \|s\| + \|t\|$ therefore condition \eqref{eq:weird_condition} is equivalent to
	\BEQ
	 \inf_w\{ \|w\| : Pw=s, \bar Pw=t\} \geq \|s\| + \|t\|,
	\EEQ
	for any $s,t \in \calE$, which reads
	\BEQ
	\|w\| \geq \|Pw\| + \|\bar P w\|,
	\EEQ
	for any $w\in \calE$.
\end{proof}

%

\end{document}